\documentclass[a4paper,9pt]{extarticle}
\usepackage[T1]{fontenc}    
\usepackage[utf8]{inputenc} 

\usepackage{lmodern}        
\usepackage{amsmath}
\usepackage{graphics}
\usepackage{latexsym}
\usepackage{amsfonts}
\usepackage{bbold}
\usepackage{amsthm}
\usepackage{amssymb}
\usepackage{stmaryrd}
\usepackage{algorithm}
\usepackage{graphicx}
\usepackage{algpseudocode}
\usepackage{color}
\usepackage{indentfirst}
\usepackage{caption}
\usepackage{subcaption}
\usepackage{titlesec}
\usepackage{cite}
\usepackage{multirow}
\usepackage{enumerate}
\usepackage[colorlinks=true, allcolors=blue]{hyperref}
\usepackage{lmodern}
\titleformat{\section}
  {\normalfont\bfseries}{\thesection}{0.5em}{}
\titleformat{\subsection}
  {\normalfont\itshape}{\thesubsection}{0.5em}{}
\numberwithin{equation}{section} 
\leftmargin
-1.25cm 
\title{
 Numerical approximations of the value of zero-sum stochastic differential impulse controls game in finite horizon 
 }
\author{Brahim EL ASRI \thanks{Universit\'e Ibn Zohr, Equipe. Aide \'a la decision,
ENSA, B.P.  1136, Agadir, Maroc. e-mail: b.elasri@uiz.ac.ma }\,\,\,
\, and \, Antoine ZOLOME \thanks{Universit\'e Ibn Zohr, Equipe. Aide \'a la decision,
ENSA, B.P.  1136, Agadir, Maroc. e-mail: kossiantoine.zolome@edu.uiz.ac.ma.} }
\begin{document}
\date{}
\maketitle
\newtheorem{theo}{Theorem}[section]
\newtheorem{problem}{Problem}[section]
\newtheorem{pro}{Proposition}[section]
\newtheorem{cor}{Corollary}[section]
\newtheorem{axiom}{Definition}[section]
\newtheorem{rem}{Remark}[section]
\newtheorem{lem}{Lemma}[section]
\newtheorem{ex}{Example}[section]
\newtheorem{Heuristics}{Heuristics}
\newtheorem{ass}{Assumption}[section]

\newcommand{\bass}{\begin{ass}}
\newcommand{\eass}{\end{ass}}
\newcommand{\brm}{\begin{rem}}
\newcommand{\erm}{\end{rem}}
\newcommand{\bethe}{\begin{theo}}
\newcommand{\eethe}{\end{theo}}
\newcommand{\bl}{\begin{lem}}
\newcommand{\el}{\end{lem}}
\newcommand{\bp}{\begin{pro}}
\newcommand{\ep}{\end{pro}}
\newcommand{\bcor}{\begin{cor}}
\newcommand{\ecor}{\end{cor}}
\newcommand{\be}{\begin{equation}}
\newcommand{\ee}{\end{equation}}
\newcommand{\beq}{\begin{eqnarray*}}
\newcommand{\eeq}{\end{eqnarray*}}
\newcommand{\beqa}{\begin{eqnarray}}
\newcommand{\eeqa}{\end{eqnarray}}
\newcommand{\dg}{\displaystyle \delta}
\newcommand{\cm}{\cal M}
\newcommand{\cF}{{\cal F}}
\newcommand{\cR}{{\cal R}}
\newcommand{\bF}{{\bf F}}
\newcommand{\tg}{\displaystyle \theta}
\newcommand{\w}{\displaystyle \omega}
\newcommand{\W}{\displaystyle \Omega}
\newcommand{\vp}{\displaystyle \varphi}
\newcommand{\ig}[2]{\displaystyle \int_{#1}^{#2}}
\newcommand{\integ}[2]{\displaystyle \int_{#1}^{#2}}
\newcommand{\produit}[2]{\displaystyle \prod_{#1}^{#2}}
\newcommand{\somme}[2]{\displaystyle \sum_{#1}^{#2}}
 \newcommand\omicron{\textit{o}}
\newlength{\inter}
\setlength{\inter}{\baselineskip}
\setlength\parindent{24pt}
\setlength{\baselineskip}{7mm}
\newcommand{\no}{\noindent}
\newcommand{\rw}{\rightarrow}
\def \ind{1\!\!1}
\def \R{I\!\!R}
\def \N{I\!\!N}
\def \cadlag {{c\`adl\`ag}~}
\def \esssup {\mbox{ess sup}}
\hrule
\begin{abstract} 
In this paper, we consider a differential stochastic zero-sum game in which two players  intervene by adopting impulse controls in a finite time horizon. We provide a numerical solution as an approximation of the value function, which turns out to be the same for both players. While one seeks to maximize the value function, the other seeks to minimize it. Thus we find a single numerical solution for the Nash equilibrium as well as the optimal impulse controls strategy pair for both player based on the classical Policy Iteration (PI) algorithm. Then, we perform a rigorous convergence analysis on the approximation scheme where we prove that it converges to its corresponding viscosity solution as the discretization step approaches zero, and under certain conditions. We showcase our algorithm by implementing a two-player 
almost analytically solvable game in which the players act through impulse control and compete over the exchange rate.
\end{abstract}
\hrule
\section{Introduction}
The numerical solution of optimal control problems has become a serious necessity in financial problems and stochastic processes theory for several applications in various fields such in the banking sector, insurance and reinsurance, mathematical finance, and management science.
Stochastic differential games (SDG) in turn constitute an interesting branch of mathematics as they allow to model the interaction between two or more agents whose objective functions depend on the evolution of a continuous time stochastic process. The foremost recent studies are centered on the stochastic impulse games (SIG) which combine both stochastic differential games and stochastic impulse controls. The first theories on the topic were pioneered by Isaacs \cite{[ISAAC]}  followed by Evans \& Souganidis \cite{[EVSOU]} and Fleming \& Souganidis (see also Lions \& Souganidis \cite{[LISOU1],[LISOU2]}). Later on, Cosso \cite{[COSSO]} and El Asri \&  Mazid \cite{[BESM]} (see also  El Asri \& al. \cite{[BELASM],[BELA]}) studied  the problem involving impulse controls with double obstacle quasi-variational inequalities where it is shown that the two (upper and lower) value functions are viscosity solutions to the associated variational equation and therefore proving that  the game admits a solution (Nash equilibrium). As a result, the solution to the stochastic control problem can be computed by solving the corresponding Hamilton-Jacobi-Bellman Inequality (HJBI) as long as the convergence to the viscosity solution is guaranteed. However all these references mentioned above were focused on the dynamic programming principle (DPP) approach and viscosity solutions frameworks.\\
\textit {1.1. Background for the study.} In our paper we focus on solving numerically this following discrete Hamilton-Jacobi-Bellman-type obstacle problem in a finite horizon:
\begin{align} \label{introHJB}
\begin{split}
0 = F(t,& x,V(t,x), DV(t,x),D^2V(t,x),H^c_{sup} V(t,x),H^\chi_{inf} V(t,x)) \\ 
&:=\left\{
\begin{array}{l}
\max\Big\{ \min\Big[-\cfrac{\partial V}{\partial t}  -\mathcal{L}V-f,V-H^c_{sup} V\Big],V-H^\chi_{inf} V \Big\}=0
\qquad \text{on  } [0,T)\times\mathbb{R}^d, \\ \\
V(T,x)=g(x)\qquad\qquad\qquad\qquad\qquad\qquad\qquad\qquad\qquad\qquad  \forall x\in\mathbb{R}^d,
\end{array}
\right.
\end{split}
\end{align}
where $\mathcal{L} = \mathcal{L}(t,x)$ is the (possibly degenerate) generator of a SDE defined by
$$\mathcal{L}V(t,x)= \langle b(t,x),\nabla _x V(t,x)\rangle +\cfrac{1}{2}\text{ trace}[\sigma\sigma^*(t,x)\nabla^2_xV(t,x)].$$  
The function $f = f(t,x)$ is a forcing term and $H^c_{sup}/H^\chi_{inf}$ are the impulses (a.k.a intervention)  operators
$$H^c_{sup}V(t,x)=\sup\limits_{\xi\in \mathcal{U}}[V(t,x+\xi)-c(t,\xi)],\qquad H^\chi_{inf}V(t,x)=\inf\limits_{\eta\in \mathcal{V}}[V(t,x+\eta)+\chi(t,\eta)].$$ 
Problem (\ref{introHJB})  describes the solution of a class of zero-sum stochastic differential games (ZSSDG) involving impulse control, where the value function is both time and space dependant. This type of problem comes up very frequently in the literature on optimal impulse control problems under various aspects. In particular, Azimzadeh \cite{[AZ]} studied SDGs with asymmetric setting in which one player  plays an impulse control while their opponent plays only a stochastic control, and with no restrictions on cost functions (see also Zhang \cite{[FZ]}, Yong \cite{[JM]}, Souquière \cite{[ASO]}, Baňas \& al. \cite{[LBGF]}). In \cite{[BELASM]} the authors considered a deterministic impulse control game in infinite time horizon, the authors in \cite{[COSSO],[BESM]}  analyse the instance when both players adopt impulse controls in respectively infinite and finite time horizon under rather weak assumptions on the cost functions. Even the more general and versatile case of  NonZero-Sum Stochastic impulses games (NZSSIG) has recently received attention (see the works of Zabaljauregui \cite{[DIEGOZ]},  Zabaljauregui \& al \cite{[DZ]},  Aïd \& al. \cite{[RAMB]}, De Santis \cite{[DDS]}). \\ 
For this type of problems, the HJBI equation  turns out to be a \textit{min-max} double-obstacle quasi-variational inequality, with the two obstacles  implicitly given and depending on the solution.
To be precise, with comparison  to a conventional  single player obstacle problem, the presence of the second player makes matter more complex and challenging to solve in the sense that
\begin{itemize}
    \item  Both players use impulse controls.\footnote{When both players use impulse controls, they must take into consideration how their interventions interact with each other and with the system. This is difficult because these impulses can occur at any time and have significant impact on the state of the system.}
    \item   The value function of the ZSSDG is the unique viscosity solution of a HJBI double obstacle problem. \footnote{These type of problems are notoriously difficult to solve, even in the case where only one player uses impulse controls.}.
    
\end{itemize}
The later implies that the value functions of the game cannot  be defined by a simple maximisation, instead the optimality  is characterized as a pair of value functions displaying the best expected outcome of each  player. The pair of strategies at which these values functions are  obtained is then known as a  \textit{Nash equilibrium}. However, this is less accentuated for the  ZSSDGs where the gains (resp, losses) undergone by one player are converted into  losses (resp, gains) for his opponent. Thus both value functions add-up to zero and the problem is reduced to  finding  one value function.  \\
 In order to compute the value function, numerical schemes have been derived and studied for a long time, we  refer for instance to Barles \& Souganidis \cite{[BS91]}, Bokanowski \& al. \cite{[BFS16]}, Camili \& Zidani \cite{[ABFCHZ]}, Barles \& Jakobsen \cite{[GBERJ]} and  also the recent works  of Azimzadeh \cite{[AZ], [AZh], [AZhBEM], [AZhBEM2]}  for the derivation of such schemes   and   for the study of their properties, including some proofs of convergence and of rate of convergence. To summarize,  we make an  attempt  at numerically solving a zero-sum stochastic impulse games where the value function is both time and space dependant and where the time horizon is finite. Our results echo those of Barles and Souganidis in the sense that to ensure  the theoretical convergence to the solution of (\ref{introHJB}), the  schemes are required to be \textit{monotone}, \textit{stable}, and satisfy \textit{consistency} condition.
 \\ 
 \no \textit{1.2. Contribution } Our paper contributes to both the theoretical understanding and practical applications of differential games with impulse controls in the the following  aspects. First, by  developing an approximation method of the HJBI equation and computational algorithms which leads to an efficient and accurate solution. Secondly, by providing a detailed analysis of the convergence properties and accuracy of the finite difference methods for solving the game. And lastly, Applying the numerical methods to real-world problems in finance and evaluating its its effectiveness in solving such problems.

We  consider for the purpose approximation schemes of the type 
\begin{equation}\label{schIntro}
    S \Big(h,(t,x),V^h , [\mathcal{H}^{h}_{sup} V^h](t,x),[\mathcal{H}^{h}_{inf} V^h](t,x)\Big) = 0,
\end{equation}
where $S$ is a consistent, monotone and  uniformly continuous approximation scheme of $F$ in (\ref{introHJB}) and the coupling among both equations is only in the variable $V^h$ where $V^h$ is expected to be an approximation to the solution $V$ of system (\ref{introHJB}). Typical approximation of this type, like ﬁnite diﬀerence methods and semi-Lagrangian schemes, will be discussed in detail. As far as we know, the results in this paper on numerical approximation  of 
the solution of the ZSSDIG  are the first of their kind. We are also aware of the work by El Asri \& Lalioui \cite{[BELA2]}  using numerical methods to study a related class of problem, namely, a deterministic finite-time horizon, two-player, zero-
sum differential games (DGs) where the \textit{maximizing} player is allowed to take \textit{continuous} \textit{and} \textit{impulse} \textit{controls} whereas the \textit{minimizing} player is allowed to take
\textit{impulse control only}. \\
\no \textit{1.3. Organization of the paper. }  

In the next two subsections, we establish the problem formulation by appropriately defining the ZSSDIG and recalling the primary result from \cite{[BESM]}, which guarantees the existence and uniqueness of a viscosity solution to the HJBQVI associated with this problem. We present the general assumptions and hypotheses that underlie the proofs of these results, for the readers to grasp the key concepts as a foundation for our subsequent analysis.  Furthermore, we demonstrate an example of application used for our numerical experiment, and conclude this section by providing a convergence analysis for our scheme, ensuring its convergence to the viscosity solution. This analysis provides valuable insights into the efficacy and reliability of our proposed approach. Section \ref{section4} introduces a discretization schemes, namely the "\textit{Semi-Lagrange explicit-scheme}". We highlight the complexity of the scheme and  analyze the respective properties, specifically the necessary monotonicity, stability, and consistency properties required for convergence to hold. By ensuring these properties, we guarantee the accuracy and reliability of the discretization process and validate the analysis provided in section \ref{section2}. In Section \ref{section3}, we delve into the policy iteration algorithm (\textit{algorithm} \ref{alg1}), discussing relevant properties and assumptions crucial to ensure the convergence of the algorithm towards an optimal solution. Finally, in Section \ref{section5}, we present  numerical evidence that serves as an illustration of  our work in this paper, in the light  of which we draw some conclusions in Section \ref{section6}.

\section{Assumptions and Formulation of the problem}\label{section2}
\no
\subsection{Two player zero-sum stochastic impulse game}\label{section21}

Since our main focus in this paper is to  provide a numerical approximation to the solution of  problem (\ref{introHJB}), we only give in this section  a  brief  description of the game   and  refer the interested reader  to \cite{[BESM]} (and the references therein) for a full detailed explanation of the settings and more in-depth information relative to  the general modelling framework  and proofs to be considered.\\
We denote by $(\Omega,\mathcal{F},(\mathcal{F}_t)_{t\leq T},\mathbb{P})$  a fixed and  filtered probability space under the usual conditions on which we define a standard $d$-dimensional Brownian motion $W = (W_t)_{t\leq T}$. We consider the two player  zero-sum game in which each player can intervene when convenient to them through a set of control actions at some non-decreasing stopping time frequencies $\{\tau_m\}$ (resp, $\{\rho_k\}$) of $[0,T]$ for player I (resp, for player II), satisfying $\tau_m \leq \tau_{m+1}$ (resp, $\rho_l \leq \rho_{l+1}$). 
We are given two convex cone sets  $\mathcal{U}$ and $\mathcal{V}$ of $\mathbb{R}^d$, with $\mathcal{V}\subset \mathcal{U}$ representing the spaces of control actions. The impulse control actions for player I (resp, for player II) take the form $u=\sum\limits_{m\geq1}\xi_{m}\ind_{[\tau_{m},T]}$ (resp,  $v=\sum\limits_{l\geq1}\eta_{l}\ind_{[\rho_{l},T]}$) on  $[t,T]\subset\mathbb{R}^{+}=[0,+\infty)$ such that $u \in  \mathcal{U}$ and $v \in \mathcal{V}$ and $\xi_m$ (resp, $\eta_l$) the actions are $\mathcal{F}_{\tau_{m}}$-measurable (resp., $\mathcal{F}_{\rho_{l}}$-measurable).\\
We assume the game starts at   time $t \in [0,T]$ (in a finite time horizon, e.g $T< +\infty$), the state variable process $X = X^{t,x,u,v}$ then with initial value $ X_{t^-} = x$ ,  $x \in \mathbb{R}^d $ viewed at a time $s$ $( t\leq s \leq T)$ involving the impulses $u$ and $v$ is described as the follow : 
\begin{equation}
\label{az}
\begin{array}{ll}
 X_{s}=x+\integ{t}{s}b(r,X_{r})dr+\integ{t}{s}\sigma(r,X_{r})dW_{r}
+\sum\limits_{m\geq 1}\xi_{m}\ind_{[\tau_{m},T]}(s)\prod_{l\geq 1}\ind_{\{\tau_{m}\ne\rho_{l}\}}+\sum\limits_{l\geq 1}\eta_{l}\ind_{[\rho_{l},T]}(s),
\end{array}
\end{equation}
where $b$ and $\sigma$ are respectively the drift and volatility of the process, and $\xi_m$, $\eta_l$ the controls actions induced  by each player.
In other words, if there is no intervention from both players, the state variable evolves independently through an Ito's diffusion process. Both  can then, at any time, attempt to shift the  process by  performing actions (the impulse controls) whenever convenient. 
The impulses actions $u$ and $v$ can then be defined as a strategy $(\tau_m,\xi_m)_{m\geq 1}$ (resp, $(\rho_l,\eta_l)_{l\geq 1}$) for player I (resp, player II). We will assume that player II has the priority when both player act simultaneously. This former statements is purely a matter of convention and is not very restrictive\footnote{If we give priority to the player I instead, the product $\prod_{l\geq 1}\ind_{\{\tau_{m}\ne\rho_{l}\}}$ will not restrict anymore the interventions of player I but will be on the side of player II control actions.
}. 
Considering a pair of strategy $ \big \{\psi := (u := (\tau_m,\xi_m)_{m\geq 1} ,v :=(\rho_l,\eta_l)_{l\geq 1}) \big \}$ and our state process dynamics, the expected payoff $ J$ (e.g, gain functional to maximize over all $\psi$  for player I and cost functional to minimize over all $\psi$ for player II) is defined as : 
\begin{equation}
\label{aar}
\begin{array}{ll}
J(t, x; \psi)=\mathbb{E}\bigg[\integ{t}{T}f(s,X_{s}^{t,x,\psi})\exp(-\lambda(s-t))ds&-
\sum\limits_{m\geq 1} c(\tau_{m},\xi_{m})\exp(-\lambda(\tau_{m}-t)) \ind_{[\tau_{m}\leq T]}\prod_{l\geq 1}\ind_{\{\tau_{m}\ne\rho_{l}\}}  \\&+  \sum\limits_{l\geq 1} \chi(\rho_{l},\eta_{l})\exp(-\lambda(\rho_{l}-t)) \ind_{[\rho_{l}\leq T]} 
 +g(X_{T}^{t,x,\psi})\exp(-\lambda(T-t))\bigg],
\end{array}
\end{equation}

where : 
\begin{itemize}
 \item $f=f(t,x)$, $(t,x) \in [0,T]\times \mathbb{R}^d$ represents the running reward; 
 \item $g$ is  the reward receive at time $T$, modelled by $g = g(x)$, $x\in \mathbb{R}^d$;
 \item $c = c(t,\xi) $ refers to the intervention cost function for player I, using control $\xi$ and analogously
 \item $\chi = \chi(t,\eta)$, the intervention cost function for player II when using control $\chi$.\\
 Note that the player performing his actions has to pay an intervention cost yielding to a gain for the other player and vice-versa (which is the essence of the \textit{zero-sum} game).
\end{itemize}
 we consider  the following assumptions 
about the various quantities appearing in our problem: 
\begin{ass}\label{ass1}
$[\textbf{H}_{b,\sigma}]$  $b:[0,T]\times \mathbb{R}^{d}\rightarrow
\mathbb{R}^{d}$ and $\sigma :[0,T]\times \mathbb{R}^{d}\rightarrow
\mathbb{R}^{d\times d}$ are two Lipschitz continuous functions for which there
exists a constant $C>0$ such that for any $t\in \lbrack 0,T]$ and
$x,x^{\prime }\in \mathbb{R}^{d}$
\begin{equation}
|\sigma (t,x)-\sigma (t,x^{\prime })|+|b(t,x)-b(t,x^{\prime })|\leq
C|x-x^{\prime }|.
 \label{eqs}
\end{equation}%
 Also there exists a constant $C>0$ such that for any $(t,x)\in[0,T]\times\mathbb{R}^d$
\begin{equation}
|\sigma (t,x)|+|b(t,x)|\leq C(1+|x|).
\end{equation}
$[\textbf{H}_{f}]$  $f:[0,T]\times\mathbb{R}^{d} \rightarrow \mathbb{R}$
is uniformly continuous and bounded on $[0,T]\times\mathbb{R}^{d}$ and the function $ g:\mathbb{R}^{d} \rightarrow\mathbb{R}$ is uniformly continuous and bounded on $\mathbb{R}^n$.\\
$[\textbf{H}_{c,\chi}]$ The cost functions $c:[0,T]\times \mathcal{U} \rightarrow \mathbb{R}$ and $\chi:[0,T]\times \mathcal{V} \rightarrow \mathbb{R}$ are measurable
and uniformly continuous. Furthermore,
\begin{equation}\label{cout} \underset{[0,T]\times \mathcal{U}}{inf}c \geq k, \hspace*{2cm} \underset{[0,T]\times \mathcal{V}}{inf}\chi\geq k,
\end{equation}
where $k>0$. Moreover,
\begin{equation} \label{cout1}
  c(t,\xi_1+\xi_2)\leq c(t,\xi_1)+c(t,\xi_2),\\
\end{equation}
\begin{equation}\label{cout2}
  \chi(t,\eta_1+\eta_2)\leq \chi(t,\eta_1)+\chi(t,\eta_2),
\end{equation}
for every $t\in [0,T], \xi_1,\xi_2\in \mathcal{U}$ and $\eta_1,\eta_2\in\mathcal{V}$.
\\
$[\textbf{H}_{g}]$ (no terminal impulse). For any $x\in\mathbb{R}^d, \eta \in\mathcal{V} $ and $\xi\in\mathcal{U}$,
\begin{equation}\label{Terminal}
\sup\limits_{\xi\in \mathcal{U}}[g(x+\xi)-c(T,\xi)]\leq g(x) \leq \inf\limits_{\eta\in \mathcal{V}}[g(x+\eta)+\chi(T,\eta)].
\end{equation}
$[\textbf{H}_{\tau,\rho}]$ The actions times $\tau_m$ and $\rho_l$ take value in a set $\mathbb{Q}_{[t,T]} =(\mathbb{Q} \cap [t,T])\cup \{t,T\}. $
\\
$[\textbf{H}_{c,\chi,\xi,\eta}]$ There exists a function $h :[0,T] \to (0,\infty)$ such that for all $t \in [0,T]$ :
\begin{equation}\label{hcxen}
\chi(t,\eta_1 +\eta_2) \leq  \chi(t,\eta_1)+\chi(t,\eta_2) - h(t),
\end{equation}
and 
\begin{equation}
    c(t,\xi_1 + \eta + \xi_2 ) \leq  c(t,\xi_1) - \chi(t,\eta) +c(t,\xi_2) -h(t),
\end{equation}
for every  $t\in [0,T], \xi_1,\xi_2\in \mathcal{U}$ and $\eta_1,\eta, \eta_2\in\mathcal{V}$.
\end{ass}
Throughout  the paper, only admissible impulse controls  were considered, which means, the strategy $\psi=(u,v)$  is admissible if it gives a well defined expectation of the payoffs for both players (i.e: a solution of (\ref{aar})  exists and is unique) and the average number of impulses ($\mu_{t,T}(u):=\sum\limits_{m\geq1}\ind_{[\tau_{m},T]} $ and $\mu_{t,T}(v):=\sum\limits_{l\geq1}\ind_{[\rho_{l},T]}$) is finite (i.e: $ \mathbb{E}[\mu_{t,T}(u)]<\infty$  and  $\mathbb{E}[\mu_{t,T}(v)]<\infty$). The set of all admissible impulse controls for player I
(resp., II) on $[t,T]$ is denoted by $\mathcal{U}_{t,T}$ (resp., $\mathcal{V}_{t,T})$. In this context, an \textit{optimal}  strategy $\psi^* =(u^*,v^*)$  is  is defined  if it  is indeed admissible and if : 
$$ J(t,x , u^*,v^*) \geq J(t,x , u,v^*) \qquad \text{and} \qquad J(t,x , u^*,v^*) \geq J(t,x , u^*,v). $$
If this  strategy exists, we define a solution of the game (\textit{value-function})  by : 
$$ V(t,x) = J(t,x , u^*,v^*) =J(t,x , \psi^*). $$ 

\begin{axiom} (nonanticipative strategy) The nonanticipative strategy set $\mathcal{A}_{t,T}$ for player I
is the collection of all nonanticipative maps $\alpha$ from  $\mathcal{V}_{t,T}$ to $\mathcal{U}_{t,T}$, i.e  for any $[t,T]$-valued $\mathbb{F}$-stopping times $\tau$ and any $v_{1}$,$v_{2}\in \mathcal{V}_{t,T}$,
\begin{equation*}
\begin{array}{ll}
\text{if}\quad v_{1}\equiv v_{2}\quad \text{on} \quad [t,\tau ]
,\quad \text{then}\\ \\
 \;\alpha(v_{1})\equiv \alpha(v_{2})\quad \text{on}\quad [t,\tau ].
\end{array}
\end{equation*}
Analogously, the nonanticipative strategy set $\mathcal{B}_{t,T}$ for player II
is the collection of all nonanticipative maps $\beta$ from  $\mathcal{U}_{t,T}$ to $\mathcal{V}_{t,T}$.
\end{axiom}
\subsection{The system of HJBI quasi-variational inequalities}
\no
Given the above assumptions, we  want to know whether the game admits a solution, and if it exists how to compute it. We are clearly aware  that one player (player I) seeks to maximize his payoffs while the opponent (player II) attempts to minimize his, this implies the definition of an upper $V^{+}(.)$ and lower $V^{-}(.)$ value functions for the game, for every $(t,x)\in[0,T]\times\mathbb{R}^d$ such that:
\begin{equation}
V^{-}(t,x):=\inf\limits_{\beta\in {\mathcal{B}_{t,T}}}\sup\limits_{u\in\mathcal{U}_{t,T}}J(t,x,u,\beta(u)),
\end{equation}
\begin{equation}
V^{+}(t,x):=\sup\limits_{\alpha\in \mathcal{A}_{t,T}}\inf\limits_{v\in\mathcal{V}_{t,T}}J(t,x,\alpha(v),v).
\end{equation}
The game is said to admit a solution if the upper and lower value functions exists and  coincide: $ V:= V^{+} =V^{-}$, we refer to  $V$  as  the value function of the game. \\
The  Hamilton-Jacobi-Bellman-Isaacs equation  associated to our problem  which turns out to be the same for the two value functions, as the two players can not act simultaneously on the system, is given by: 

\begin{equation} \label{eq:HJBI}\left\{
\begin{array}{l}
\max\Big\{ \min\Big[-\cfrac{\partial V}{\partial t} -\mathcal{L}V-f,V-H^c_{sup} V\Big],V-H^\chi_{inf} V \Big\}=0
\qquad \text{on  } [0,T)\times\mathbb{R}^d, \\ \\
V(T,x)=g(x)\qquad\qquad\qquad\qquad\qquad\qquad\qquad\qquad\qquad\qquad  \forall x\in\mathbb{R}^d.
\end{array}
\right.
\end{equation}
where $\mathcal{L}$, the infinitesimal generator of $X$ when uncontrolled, is the second-order local operator
   $$\mathcal{L}V(t,x)= \langle b(t,x),\nabla _x V(t,x)\rangle +\cfrac{1}{2}\text{ trace}[\sigma\sigma^*(t,x)\nabla^2_xV(t,x)],$$ 
and the non-local intervention operators  $H_{sup}^c$ (the cost operator) and $H_{inf}^\chi$ (the gain operator) are given by : 
$$H^c_{sup}V(t,x)=\sup\limits_{\xi\in \mathcal{U}}[V(t,x+\xi)-c(t,\xi)],\hspace*{1cm} H^\chi_{inf}V(t,x)=\inf\limits_{\eta\in \mathcal{V}}[V(t,x+\eta)+\chi(t,\eta)],$$
for every $(t,x)\in[0,T]\times\mathbb{R}^d$ and 
$V$ continuous in time and space.
\begin{rem}
It is shown in \cite{[BESM]} that the equation (\ref{eq:HJBI})  admits a solution in viscosity sense that is unique in the space of bounded continuous functions on $[0,T]\times R^d$. All the proofs on the existence and uniqueness of the viscosity solution as well as additional assumptions and property of the value function is given in that paper.
\end{rem}
 We recall the formulations below for the corresponding  notions of viscosity solution and some useful properties related to the sequel of our paper.
\begin{axiom}
We say  $V =V(t,x)$ is a viscosity subsolution of (\ref{eq:HJBI}) if $V$ is upper semi-continuous and if  for any $(t_0,x_0)\in[0,T)\times\mathbb{R}^d$ and  any function
 $\phi \in C^{1,2}([0,T)\times\mathbb{R}^d)$, such that $(t_0,x_0)$ is a local maximum  of $V-\phi$, we have:
\begin{equation}
\begin{array}{ll}
\max\Big\{ \min\Big[-\cfrac{\partial \phi}{\partial t}(t_0,x_0) + \lambda \phi(t_0,x_0)-\mathcal{L}\phi(t_0,x_0)-f(t_0,x_0),V(t_0,x_0)-H^c_{sup} V(t_0,x_0)\Big],\\\qquad\qquad\qquad V(t_0,x_0)-H^\chi_{inf} V(t_0,x_0) \Big\}\leq 0.
\end{array}
\end{equation}

We say $V =V(t,x)$  is  a viscosity supersolution of (\ref{eq:HJBI}) if $V$  is lower semi-continuous, and  if for any $(t_0,x_0)\in[0,T)\times\mathbb{R}^d$ and any function
 $\phi \in C^{1,2}([0,T)\times\mathbb{R}^d)$, such that  is  $(t_0,x_0)$  a local
minimum of $V-\phi$, we have:
\begin{equation}
\begin{array}{ll}
\max\Big\{ \min\Big[-\cfrac{\partial \phi}{\partial t}(t_0,x_0) + \lambda \phi(t_0,x_0)-\mathcal{L}\phi(t_0,x_0)-f(t_0,x_0),V(t_0,x_0)-H^c_{sup} V(t_0,x_0)\Big],\\\qquad\qquad\qquad V(t_0,x_0)-H^\chi_{inf} V(t_0,x_0) \Big\}\geq 0.
\end{array}
\end{equation}

We say $V =V(t,x)$  is  a viscosity solution   if $V$ is both a viscosity supersolution and subsolution.
\end{axiom}
\no
Roughly speaking, the notion of viscosity solution is obtained by replacing the derivatives of  the function $V$ with the test function $\phi$ such that :
\begin{equation}\label{visco1}\tag{\textit{definition-1}}
     F(t, x,V(t,x), D\phi(t,x),D^2\phi(t,x),H^c_{sup} V(t,x),H^\chi_{inf} V(t,x)) = 0.
\end{equation}
Another notion among all the several ways to establish the viscosity solution is to replace instead all the function $V$ (including also the intervention operators) by the test function $\phi$ such that:
\begin{equation}\label{visco2}\tag{\textit{definition-2}}
     F(t, x, \phi(t,x), D\phi(t,x),D^2\phi(t,x),H^c_{sup} \phi(t,x),H^\chi_{inf} \phi(t,x)) = 0.
\end{equation}
\begin{rem}\label{remv12}
 These two concepts of the viscosity solution are not in general equivalent but have the same meaning when imposing further properties on the problem. Such results were developed in \cite{[ALTO]} for a class of particular nonlocal operators (integro-differential operators).  It is also shown in \cite{[GW09]} that for a uniformly continuous function, a viscosity solution in (\ref{visco2}) sense is also a  solution in (\ref{visco1}) sense. And when the interventions operators are nondecreasing, a semicontinuous solution in (\ref{visco1}) sense is also a solution in  (\ref{visco2}) sense, as a result, a comparison principle in (\ref{visco2}) sense implies a comparison principle in (\ref{visco1}) sense. However for nonlocal equations arising from impulse control, the most relevant notion of viscosity solution are established from a (\ref{visco1} ) sense (see, \cite{[AZhBEM]} for more details).
\end{rem}
In all the  following, we  assume Assumption \ref{ass1} holds. 
\begin{theo} (Comparison theorem) The HJBQVI (\ref{eq:HJBI}) satisfies a comparison principle in the set of  bounded maps  valued from $[0,T]\times \mathbb{R^n}$ (for any $ n\geq 1$) (i.e., if $u$ is a bounded uniformly continuous viscosity subsolution of (\ref{eq:HJBI}) and $w$ is a  bounded uniformly continuous viscosity supersolution of (\ref{eq:HJBI}) ,then $u \leq w$).
\begin{proof}
    Established in  \cite{[BESM]}.
\end{proof}
\end{theo}
Unless otherwise specified, a comparison principle will be define in (\ref{visco1}) sense.
\begin{pro} (\cite{[BESM]}, Proposition 3) \label{propopo}
The upper and lower value functions satisfy the following properties for all $(t,x) \in [0,T)\times \mathbb{R}^d$.
\begin{itemize}
    \item[(i)] $H_{inf}^\chi V(t,x) \geq V(t,x)$,
    \item[(ii)] If $H_{inf}^\chi V(t,x) > V(t,x)$ then $H_{sup}^c V(t,x) \leq V(t,x)$.
\end{itemize}
\end{pro}

\begin{cor} (\cite{[BESM]}, Corollary 3)
The lower and upper value functions  are continuous and coincide, and the value function  of the stochastic diﬀerential game,unique among all the bounded solutions is given by $V(t,x) := V^{-}(t,x) = V^{+}(t,x)$ for every $(t,x) \in[0,T) \times \mathbb{R}^d$.
\end{cor}
\subsection{Example of application:} 
The exchange rate serves as a representation of the value of a foreign country's currency in relation to the domestic country's currency.  In theory, the management of the exchange rate is a significant tool for the Central Bank in its effort to reduce high fluctuations in the exchange rate. By maintaining the exchange rate as close as possible to a predetermined range or target, the Central Bank can regulate its volatility. \\
The numerical application  in this research  considers the exchange rate problem  in which a commercial institution (P1) and a Central Bank (P2) intervene in the foreign exchange (forex) market.\\
{\Large Scenario Description }
\
\begin{enumerate}
    \item Players P1 and P2

\begin{itemize}
    \item The commercial institution (P1) aims to maximize profits by exploiting short-term fluctuations in the exchange rate between the two currencies (e.g., USD/EUR) . He can buy or sell the foreign currency (impulse control).
    \item The Central Bank (P2) manages the exchange rate to maintain stability within a certain range. They can intervene in the market by buying or selling their own currency (impulse control) to influence the exchange rate.
\end{itemize}
\item Objectives functions 
\begin{itemize}
    \item The commercial institution (P1): Maximizes the profits based of the rate movements.
    \item The Central Bank (P2):  Minimizes the deviation of the rate from a target value. 
\end{itemize}
\item Impulses controls \\
     Both players can make discrete interventions in the market at specific times to influence the exchange rate. The intervention of the  P1 (resp,  P2), is a double sequence of stopping times $(\tau_m)_{m\geq 0}$ (resp, $(\rho_l)_{l\geq 0}$ representing the intervention times  and  controls $(\xi_m)_{m\geq 0}$ (resp, $\eta_l)_{l\geq 0}$ representing the impulses, that is a purchase ($\xi_m,\eta_l \geq 0$) or a  sale ($\xi_m,\eta_l \leq 0$) of the other currency.

\item Game Setup 
\begin{itemize}
    \item  Time horizon T, 
    \item  Known initial value $X_0$, 
    \item Terminal reward $g(X_T)$.
\end{itemize}
\item Exchange Rate Dynamics \\
Let $ X(t)_{t\geq 0} $ be  the market dynamic of the rate of change of the currencies at time $t$. When $X_t$ is high then the foreign currency is weak and when $X_t$ is low then the currency is strong. In the absence of impulse, we assume  $X(t)$  is a stochastic process of the form: 
\begin{align}\label{exchange}
\begin{split}
     dX(t) =  -u(t)\cdot X(t^-) \cdot dt + \sigma (t) \cdot  X(t^-) \cdot dW(t) &,\\ 
                                    X(0) = x >0,
                                    \end{split}
\end{align}
where $\mu(\cdot)$ and $\sigma(\cdot)$ are both non-negative functions  representing the drift speed and volatility of the exchange rate over the time (we assume later these parameters are constants) and $W_t$  a Wiener process. In  conformity to section \ref{section21}, if both players intervenes simultaneously, only  the actions of P2 will be taken into  consideration. No restrictions are imposed as for the amounts the players can buy or sell.
\item Cost/Profit Functions 
\begin{align}\label{game}
    J^\psi(t,x) = E^{(t,x)}\Bigg[-\integ{s}{T}  e^{-\theta (s-t)}(f(X_x^\psi(t))dt -  \sum\limits_{m\geq 1} e^{-\theta (\tau_m-t)} c(\tau_{m},\xi_{m}) \ind_{[\tau_{m}\leq T]} & \prod_{l\geq 1}\ind_{\{\tau_{m}\ne\rho_{l}\}}+ \\ & \sum\limits_{l\geq 1} e^{-\theta (\rho_l-t)}\chi(\rho_{l},\eta_{l}) \ind_{[\rho_{l}\leq T]}    \Bigg],
\end{align}
\end{enumerate}
where $ f(x) = -(x-x^*)^2 $ is the running cost function/profit when the exchange rate $X$ moves away from the target exchange rate $X^*$, $\theta$ is the discount rate, $c$ and $\chi$ the intervention cost function.
\begin{equation}\label{exratecosts}
    \begin{cases}
        c(\tau_m,\xi_m)= \lambda_1|\xi_m| + k_1, \quad  \lambda_1 \geq 0; k_1 >0,\\
        \chi(\rho_l,\eta_l)=\lambda_2|\eta_l| + k_2,\quad  \lambda_2 \geq 0; k_2 >0.
    \end{cases}
\end{equation} 
Now lets suppose that the controlled exchange rate process with an initial value $ X(0) = x$ is denoted by $X_x^\psi(t)$ (where $\psi$ is  the combined impulses strategies  and is defined as follows.

\begin{equation}\label{exchangeimpulse0}
\begin{cases}
    X_x^\psi(t) = X_x(t), \quad  \; 0\leq t< \min(\tau_m,\rho_l) \\ 
    \dot{ X}_x^\psi(t) = \mu(t) ,\quad  t\neq \tau_m, t\neq \rho_l, \quad t\in [0,T]\\
     X_x^\psi(\rho_l^+) = X_x^\psi(\rho_l^-) + \eta_l, \quad l= 1,2,3, \cdots ; \\
     X_x^\psi(\tau_m^+) = X_x^\psi(\tau_m^-) + \xi_m \ind_{\{\tau_m \ne \rho_l\}}, \quad m = 1,2,3,\cdots .
     \end{cases}
\end{equation}
The Drift and volatility  are not affected  when the exchange rate process is shifted by impulse. Consequently, following a  control, the exchange rate process follows  the jump diffusion dynamics $( X_x(t))$ until one of the  players finds it necessary to intervene again.

\begin{rem}
The constants $(k_{i,i\in \{1,2\}})$ represents a fixed basic cost of each impulse control intervention  independently of the impulse itself. This implies a cost of at least $k_{i,i\in [1,2]}$ to avoid for the players the possibility to apply such interventions continuously.
\end{rem}

\subsection{Convergence Analysis}
In the sequel, the approximation scheme of the nonlocal and nonlinear PDE (\ref{eq:HJBI}) take the form 
\begin{equation}{\label{sch}}
   S \Big(h,(t,x),V^h , \mathcal{H}^{h}_{sup} V^h,\mathcal{H}^{h}_{inf} V^h\Big) = 0 \quad \text{for } \quad (t,x) \in [0,T]\times \overline{\Omega},
\end{equation} 
where $\overline{\Omega}$ is a subset of $\mathbb{R}^d$ ( $\overline{\Omega} \subset \mathbb{R}^d $ ). The function $S$  is a real-valued function mapping from $(0,\infty)\times[0,T]\times \overline{\Omega} \times B(\overline{\Omega})\times \mathbb{R}\times \mathbb{R}$ that refers to the approximation scheme.    $B(\overline{\Omega})$  is the set of bounded real-valued maps from $\overline{\Omega}$ and the 
operators $\mathcal{H}^{h}_{sup}$ and $\mathcal{H}^{h}_{inf}$ maps $B(\overline{\Omega})$ to some subset of itself and will serve here  as the approximations of the intervention operators. Intuitively $V^{h} \in B(\overline{\Omega})$ is a  solution of the scheme if it satisfies (\ref{sch}). 

It is shown that as long as a scheme is \textit{stable}, \textit{monotone} and \textit{consistent}, it converges to the unique solution  of the PDE it derives from, provided that the PDE satisfies a comparison principle in (\ref{visco2}) sense. This statement is known as the  Barles-Souganidis framework  \cite{[BS91]}. This framework provides a general approach for proving the convergence of ﬁnite difference schemes to the viscosity solution their PDE under  three sufficient conditions (\textit{monotonicity}, \textit{stability}, and \textit{consistency}). However the approach requires the stronger comparison principle in (\ref{visco2}) sense which is not developed in the literature for the HJBQVI (see Remark \ref{remv12}). In \cite{[AZh]} (see also \cite{[AZhBEM],[AZhBEM2]}), the authors discussed the issue  and establish the possibility to obtain  a convergence result for non local PDE.
They showed  that using an alternative notion of \textit{consistency} (a \textit{nonlocal} \textit{consistency}) only a comparison principle in (\ref{visco1}) is required to ensure convergence, while the \textit{monotonicity} and \textit{stability} remains identical to those in the Barles-Souganidis framework.

\begin{axiom} 
\textit{\textbf{-Monotonicity:}} a numerical scheme $S$ is monotone if it maintains the order of magnitudes in the solution. That is  for  any input values $u$ and $w$ such that $ u\geq w$ pointwise, $S(\cdot,\cdot,\cdot,u,\cdot,\cdot) \leq S(\cdot,\cdot,\cdot,w,\cdot,\cdot)$ pointwise. \\
The monotonocity condition states that the scheme is nondecreasing in time and non-increasing in $V$.  \\
\textit{\textbf{- Stability:}}
 a numerical scheme $S$ is stable if 
$
\text{there exists a unique  solution } v^h \text{ of (\ref{sch}) for each } h >0 \text{ and }$ $\sup||v^h||_{\infty} < \infty.
$\\
The stability property reads that  there exits a solution $ v^h $ with a bound independent of $h$ and that small changes in the discretization parameters do not lead to significant and unpredictable alterations in the solution over time. This property is essential for the reliability and accuracy of the numerical simulation, as it prevents numerical errors from accumulating and causing the solution to deviate from the true solution of the underlying differential equation.\\
\textit{\textbf{- Consistency:}}
 A scheme $S$ is   consistent in sense of Barles-Souganidis framework if for all $x \in \overline{\Omega}$ and $\phi \in C^2( \overline{\Omega})$ 
 
 \begin{equation}
\liminf \limits_{\substack{h\to 0 \\s,y\to t, x \\ \epsilon \to 0  }} S(h,(s,y),\phi +\epsilon,\cdot ,\cdot) \geq  F_{*}(t,x,\phi,D\phi,D^2\phi,
\cdot,\cdot),
\end{equation}

\begin{equation}
\limsup \limits_{\substack{h\to 0 \\s,y\to t, x \\ \epsilon \to 0  }} S(h,(t,y),\phi +\epsilon,\cdot ,\cdot) \leq F^{*}(t,x,\phi,D\phi,D^2\phi,
\cdot,\cdot),
\end{equation}
 where the functions $F^{*}$  and  $F_{*}$ indicate respectively  the upper  semi-continuous    and  lower
semi-continuous  envelopes of the  function $F$. 
 \begin{rem}
The consistency property in the Barles-Souganidis framework  holds only for local PDEs, however the HJBQVI (\ref{sch}) is not local due to the intervention operators. This issue is addressed in \cite{[AZh]} where the consistency is extended to a nonlocal consistency for nonlocal PDEs. Their result reads that: \\
A scheme $S$ is non-locally  consistent  if for each family $(V^h)_{h>0}$ of uniformly bounded real-valued maps from   $[0,T]\times \overline{\Omega}$ and 
for any continuous function $\phi \in C^2( \overline{\Omega})$ and $x \in \overline{\Omega}$:

\begin{equation}\label{436}
\liminf \limits_{\substack{h\to 0 \\s,y\to t, x \\ \epsilon \to 0  }} S(h,(s,y),\phi +\epsilon,\mathcal{H}^{h}_{sup} V^h,\mathcal{H}^{h}_{inf} V^h) \geq  F_{*}(t,x,\phi,D\phi,D^2\phi,
H^c_{sup} \overline{V},H^\chi_{inf} \overline{V}),
\end{equation}

\begin{equation}\label{437}
\limsup \limits_{\substack{h\to 0 \\s,y\to t, x \\ \epsilon \to 0  }} S(h,(t,y),\phi +\epsilon,\mathcal{H}^{h}_{sup} V^h,\mathcal{H}^{h}_{inf} V^h) \leq F^{*}(t,x,\phi,D\phi,D^2\phi,
H^c_{sup} \underline{V},H^\chi_{inf} \underline{V}),
\end{equation}
where $\overline{V}$ and $\underline{V}$ are  respectively the upper and lower half-relaxed limits of the family $(V^h)_{h>0}$ defined by :
$$
\overline{V}:= \limsup \limits_{\substack{h\to 0 \\ s,y \to t,x}}V^h(s,y)
\quad \text{ and } \quad \underline{V}:= \liminf \limits_{\substack{h\to 0 \\ s,y \to t,x}}V^h(s,y).
$$
The functions $F^{*}$  and  $F_{*}$    are respectively  the upper  semi-continuous    and  lower
semi-continuous  envelopes of the  function $F$. 

When there is no nonlocal operators, the nonlocal consistency is equivalent to the consistency in the Barles-Souganidis framework. This property ensures  that  sufficiently smooth solutions of the scheme will converge toward the solution of $F$.
 \end{rem}
 
\end{axiom} 
\no

\begin{theo} 
Let S be a monotone, stable, and non-locally  consistent scheme whose limiting
equation (\ref{eq:HJBI}) satisfies a  comparison principle. Then, as $h \downarrow 0$, the solution
$V^{h}$ of (\ref{sch}) converges locally uniformly to the unique  solution of  (\ref{eq:HJBI}).
\end{theo}

\begin{proof}
The above is an immediate corollary of the convergence results for  nonlocal PDEs given in   \cite{[AZh]}, \cite{[AZhBEM]}.
\end{proof}


\section{Discretization and Numerical approximation}\label{section4}
\label{discretization}

In this chapter, we introduce a numerical  approximation scheme to the solution of the HJBQVI  equation. We mainly prove that the approximation scheme for any discretization step $h$  has a unique solution $V^h$ that converges locally uniformly towards the value function as $h$ goes to zero.   We restrict our attention to the  one-dimensional case ($d=1$). 
In that case, the second order operator $\mathcal{L}V$ is reduced to : 
     $$\mathcal{L}V(t,x)=  b(t,x)\mathcal{D}_x V(t,x) +\cfrac{1}{2}\sigma(t,x)^2\mathcal{D}^2_{xx}V(t,x).$$ 
We keep in mind that since the game has a terminal conditions occurring at time $T$, the expected numerical methods should proceed backwards in time to achieve the solution at time $t=0$ as well as at any  time $s$ in  $[0,T]$. Hence, it is quite natural to refer to procedures that evolve a numerical solution from time $T$  to an earlier time $\tau^n =T -n\Delta t$ where $n$ denotes the $n^{th}$ time-step.\\
Starting with  this viewpoint, we suppose a uniformly spaced grid point $\{n\Delta t\}_{n=0}^{N} $ in time and $\{x_i\}_{i=0}^{M} $ in space with $\Delta t=h= \frac{T}{N}$, (e.g: $T =N \Delta t$) and $x_i \in [ x_0, x_M] $,  ($0\leq i \leq M$, $M > 0$) so that the approximate numerical solution  is computed on the numerical domain $[0,T]\times [ x_0, x_M]$.

\begin{rem}
For the fixed $h >0$,  the numerical schemes in  the below only define the numerical solution $V^h$ that converges to $V$ as $h$ goes to zero. Intuitively,   $h$ is the parameter  that controls the accuracy  of the schemes, a smaller value of  $h$ will correspond to a finer numerical grid and vice-versa.

\end{rem}
Then let the quantity $V_i^n = V^h(\tau^n, x_i)$ be the numerical solution at time $\tau^n = T-n\Delta t$ and on the space point $x_i$. At each time step $n$ , $V^h$ is valued in $\mathbb{R}^M$: $V^n = (V_0^n,\dots , V_M^n)$. All things being equal, we use the shorthand $b_i^n = b(\tau^n,x_i)$, $\sigma^n_i =\sigma(\tau^n,x_i)$ and $f_i^n= f(\tau^n,x_i)$ to simplify our presentation.
\begin{itemize}
    \item Time derivatives 
    \begin{equation} \label{dvt}
           \mathcal{D}_tV(\tau^n,x_i) = \cfrac{\partial V}{\partial t}(\tau^n,x_i) = \cfrac{V(\tau^n + \Delta t, x_i) -V(\tau^n, x_i) }{\Delta t} = \cfrac{V_i^{n-1} - V_i^n}{\Delta t}.
    \end{equation}
  
    \item Spatial derivatives \\
We choose a forward $\mathcal{D} _x^{+} V$ (resp, backward $\mathcal{D} _x^{-} V$) difference for  the first order derivative $\mathcal{D} _x V$   depending on the coefficient $b_i^n$ in $\mathcal{L}V$ whether it is positive (resp, negative) and such that $\cfrac{\partial V}{\partial x} =\mathcal{D} _x V = \mathcal{D} _x ^{\pm} V $
    \begin{equation} \label{dvxp}
           \mathcal{D} _x^+ V(\tau^n,x_i) = \cfrac{V(\tau^n, x_{i+1}) -V(\tau^n, x_i) }{x_{i+1} - x_i} = \cfrac{V_{i+1}^{n} - V_i^n}{\Delta x_i} \text{ in } \{b_i^n > 0\},
    \end{equation} 
      \begin{equation} \label{dvxm}
\mathcal{D} _x^- V(\tau^n,x_i) = \cfrac{V(\tau^n, x_{i}) -V(\tau^n, x_{i-1}) }{x_{i} - x_{i-1}} = \cfrac{V_{i}^{n} - V_{i-1}^n}{\Delta x_{i-1}} \text{ in } \{b_i^n \leq 0\}.
    \end{equation}
    Then : 
    \begin{equation*}
        b_i^n (\mathcal{D} _xV^n)_i = \begin{cases}
 b_i^n (\mathcal{D} _x^+V^n)_i &\text{ if  } b_i^n > 0, \\
 b_i^n (\mathcal{D}_x^-V^n)_i &\text{ if  } b_i^n \leq  0.
\end{cases}
    \end{equation*}
    And the second order derivatives $\mathcal{D}^2_{xx}V$ by : 
    \begin{equation}\label{d2vx}
        \mathcal{D}^2_{xx}V (\tau^n,x_i) =  \cfrac{V^n_{i+1} - V^n_i }{(x_{i+1} - x_i)(x_{i+1} - x_{i-1)}} - \cfrac{V^n_i - V^n_{i-1}}{(x_{i} - x_{i-1})(x_{i+1} - x_{i-1})} .
    \end{equation}
    In case of uniformly spaced grid point $\Delta x = x_{i+1} -x_i$ $ \forall i \in [0,M]$, $\mathcal{D}^2_{xx}$ is reduced to : 
    \begin{equation*} 
     \mathcal{D}^2_{xx}V (\tau^n,x_i) =  \dfrac{V_{i+1}^n - 2V_i^n + V_{i-1}^n }{(\Delta x)^2}.
    \end{equation*}
\end{itemize}


\subsection{Intervention operators}
Now, let's  discuss the discretization pattern we choose for  the  interventions operators  $H_{sup}^c$ and $H_{inf}^\chi$. The fact that the two points $(x+\xi)$ and $(x+\eta)$ might eventually be outside our grid point or  not necessarily  belongs to the numerical discrete space grid $\{x_0,\dots,x_M\}$ (e.g: $x_k<x+\xi< x_{k+1})$  leads us to  reconsider a  scheme that takes into consideration any point in the continued space $[x_0,x_M]$. Hence for both operators, a convenient discretization  will require an  interpolation method to approximate   
 $V(t, x+ \xi)$ and $V(t,x+\eta)$. This situation can be avoided when setting a discrete grid for the control spaces such that the $(x+\xi)$ and $(x+\eta)$ belongs  to the space grid $\{x_0,\dots,x_M\}$. But here we discuss the eventuality to  serve for a more general case.  We  define : 
$$
  \tilde{V}(\tau^n,x_i) = \text{Interpolation}(V,(\tau^n,x_i)). \\ \\
$$ 
Where $\text{Interpolation}(V,(\tau^n,x_i))$ is the value of the numerical solution $V^h$ at any ($\tau^n,x_i$) as approximated by a standard monotone linear interpolation : 
$$
\tilde{V}(\tau^n,x_i)= \alpha V(\tau^n,x_{k+1}) +(1-\alpha) V(\tau^n,x_k) \qquad \forall x \in [x_0,x_M], $$ 
such that:  $\qquad \exists! k / x_k\leq x <x_{k+1}$ ,$\qquad$ and $\qquad \alpha = \cfrac{x-x_k}{x_{k+1}- x_k}$.

 We want a discretization  of $H_{sup}^c$ and $H_{inf}^\chi$ that fully converges in the domain of  the control sets $\mathcal{U}$ and $\mathcal{V}$. Since these sets are possibly not finite,  we introduce the subsets $\mathcal{U}^h$ and $\mathcal{V}^h$ to serve as finite, non-empty,  approximation sets of the impulse control sets $\mathcal{U}$ and $\mathcal{V}$. We purposely define $\mathcal{U}^h$ and $\mathcal{V}^h$   to be finite to ensure that the \textit{supremum}  and \textit{infinimum} are reached and can be computed (see, proposition  \ref{Haudsorff} ) . It is understood that controls $\xi$ (resp., $\eta$) that cause $x+\xi$ (resp., $x+\eta$) to exit the numerical grid are not included in $\mathcal{U}^h$ (resp., $\mathcal{V}^h$). Therefore :
 $$
 (H^{c}V^n)_i = \sup \limits_{ (x+\xi) \in \mathcal{U}^h}[\tilde{V}(\tau^n, x_i + \xi)-c(\tau^n, \xi)],
 $$
 $$
 (H^{\chi}V^n)_i = \inf \limits_{(x+
 \eta)\in \mathcal{V}^h}[\tilde{V}(\tau^n,x+\eta) +\chi(\tau^n,\eta)].
 $$ 

Using the above notation,  $H_{inf}^\chi V(\tau^n,x_i) \approx (H^{\chi}V^n)_i$ and $H_{sup}^cV(\tau^n,x_i) \approx (H^{c}V^n)_i$. 

\subsection{Boundary and Terminal Conditions}
In order to have  a  completely well-specified  numerical space grid  problem, we need to provide boundary
conditions. A number of terminal boundary conditions can be used. Note that since we will be
solving backwards in time and  the final state occurs at $\tau^n = 0$, the only conditions to  be considered in time is $V(T,x) = g(x)$. That is, for all spatial grid points $x_i$, we set $V^0_i= g(x_i)$. 
\subsubsection{Boundary conditions}
For  spatial boundaries,  we set $\mathcal{D}^{\pm}_{x}V$ and $\mathcal{D}^2_{xx} V$  to be zero at the extreme points $x_0$ and $x_M$ as for a Newman boundary  condition reflecting  that the solution itself remain unchanged across the boundary.  However, the techniques in this paper  are not reliably specific to this choice of
boundary condition. That being said : 
$$
\mathcal{D}^{\pm}_x V(t,x_0) = \mathcal{D}^{\pm}_x V(t,x_M) = 0, \qquad \forall t \in [0,T].
$$
$$
\mathcal{D}^{2}_{xx} V(t,x_0) = \mathcal{D}^{2}_{xx} V(t,x_M) = 0,  \qquad \forall t \in [0,T].
$$
For the intervention operators, we should make sure no extrapolation is performed for  any point not contained in the grid. Then
$$
  \tilde{V}(\tau^n,x) =\begin{cases}
 V(\tau^n,x_0) &\text{ if  } x\leq x_0, \\
V(\tau^n,x_M) &\text{ if  } x \geq  x_M.
\end{cases}
$$
\subsubsection{Terminal condition}
The terminal condition $V(T,x) = g(x)$ is a direct consequence of the assumption (\ref{ass1}, $[\textbf{H}_{g}]$) which has the interpretation that performing an action at the final time $T$ does not increase the value of the payoff. This condition can be written equivalently as 
\begin{equation}\label{HgTerminal1}
    V(T,\cdot) = \min \bigg \{ \max  \big \{ H^c_{sup}V(T,\cdot), g \big \},H^\chi_{inf}V(T,\cdot) \bigg \} ,
\end{equation}
or 
\begin{equation}\label{HgTerminal2}
    \max \bigg \{ \min  \big \{  V(T,\cdot) - g, V(T,\cdot) - H^c_{sup}V(T,\cdot) \big \}, V(T,\cdot) 
 - H^\chi_{inf}V(T,\cdot) \bigg \} =0,
\end{equation}
with the assumption that $H^c_{sup}g \leq g \leq H^\chi_{inf}g$ where 
$$H^c_{sup}g(x) =\sup\limits_{\xi\in \mathcal{U}}[g(x+\xi)-c(T,\xi)]  \qquad \text{ and } \qquad H^\chi_{inf}g(x) =  \inf\limits_{\eta\in \mathcal{V}}[g(x+\eta)+\chi(T,\eta)].$$
\subsection{The numerical approximation  scheme}
We are now ready to introduce our numerical scheme, the semi-Lagrangian explicit-impulse control scheme \cite{[ZCHEN]}. First, we approximate the term $\mathcal{D}_t+ b(.)\mathcal{D}_x$ by \footnote{Here the notation $\tilde{V}^{n-1}(\cdot,x_i + hb^n_i)$ denotes the interpolation of $V^{n-1}$ at the space point $x_i + hb^n_i$. }: 
\begin{equation} \label{lagrang}
\begin{split}
    \mathcal{D}_tV(\tau^n,x_i)+ b(\tau^n,x_i)\mathcal{D}_xV(\tau^n,x_i) & \approx \cfrac{ \tilde{V}(\tau^{n-1},x_i + hb^n_i) - V(\tau^n,x_i) }{h} ,\\
                                            & \approx \cfrac{\tilde{V}^{n-1}(\cdot,x_i + hb^n_i) - V^n(\cdot,x_i)}{h}.
\end{split}
\end{equation}
The idea behind this approximation choice is to  take advantage of the fact that the first and mainly the  second derivative coefficient  are independent of the controls $(u,v)$ (i.e, $\sigma(t,x,u,v) = \sigma(t,x)$ and $b(t,x,u,v)= b(t,x)$), and  yet the time horizon $T$ is also finite,  but also allows us to get rid of the first order derivatives. This is obtained by using a  Lagrangian argument involving
path tracing of particles (\cite{[ZCHEN]}, \textit{Section 2.3.1}). 


Using  approximation (\ref{lagrang}), we define our numerical scheme $ S \Big(h,(t,x),V^h , \mathcal{H}^{h}_{sup} V^h,\mathcal{H}^{h}_{inf} V^h\Big)$ (\ref{sch}) as : 
\begin{align} \label{lagrange1st}
\begin{split}
S \Big(h&,(t,x),V^h , \mathcal{H}^{h}_{sup} V^h,\mathcal{H}^{h}_{inf} V^h\Big) = 
\\
&\left \{
\begin{array}{l}
 \max\Big\{ \min\Big[-(\tilde{V}^{n-1}(\cdot,x_i + hb^n_i) - V^n(\cdot,x_i)) -\cfrac{1}{2}h(\sigma^n_i)^2(\mathcal{D}^2_{xx}V^n)_i-hf^n_i , V^n_i-(H^{c}V^{n-1})_i 
\Big] 
                                                           , V^n_i-(H^{\chi}V^{n-1})_i 
                                                           \Big\}=0,\\ \\
 V^0_i=g(x_i).
\end{array}
\right .
\end{split}
\end{align}
Note that for the impulses controls, we  have not  used $V^n$ but
rather $V^{n-1}$ as an argument to the discretized intervention operators, hence the name \textit{explicit-impulse  }scheme. Next, an $\omicron_i^n$ term is added to the equation (see remark (\ref{remarkomicron})) such that
\begin{equation} \label{lagrange2nd}
\begin{aligned}[b]
\max\Big\{ \min\Big[-(\tilde{V}^{n-1}(\cdot,x_i + hb^n_i) - V^n(\cdot,x_i)) -\cfrac{1}{2}h(\sigma^n_i)^2(\mathcal{D}^2_{xx}V^n)_i-hf^n_i &, V^n_i-(H^{c}V^{n-1})_i -\omicron_i^n\Big] 
\\ 
                                                           &, V^n_i-(H^{\chi}V^{n-1})_i - \omicron_i^n\Big\}=0,
\end{aligned}
\end{equation}
where  $\omicron_i^n$  represents an error in approximating the intervention
operators by their  discretized intervention operator.
Then, we obtain
\begin{align} \label{lagrange11st}
\begin{split}
S \Big(h,(t,x),&V^h , \mathcal{H}^{h}_{sup} V^h,\mathcal{H}^{h}_{inf} V^h\Big) = 
\\
&\left \{
\begin{array}{l}
\begin{aligned}[b]
\max\Big\{ \min\Big[-(\tilde{V}^{n-1}(\cdot,x_i + hb^n_i) - V^n(\cdot,x_i)) -\cfrac{1}{2}h(\sigma^n_i)^2(\mathcal{D}^2_{xx}V^n)_i-hf^n_i &, V^n_i-(H^{c}V^{n-1})_i -\omicron_i^n\Big] 
\\ 
                                                           &, V^n_i-(H^{\chi}V^{n-1})_i -\omicron_i^n\Big\}=0,
\end{aligned}
\\ \\
 V^0_i=g(x_i).
\end{array}
\right .
\end{split}
\end{align}
\begin{pro}
    The approximate HJBQVI (\ref{lagrange11st}) is equivalent to the following equation:
\begin{align} \label{proposition31}
\begin{split}
\left \{
\begin{array}{l}
\begin{aligned}[b]
\min\limits_{i \in \{0,1\}} 
       \Bigg\{ (1-i) \max\limits_{j\in \{ 0,1\}} \bigg[(1-j)(\tilde{V}^{n-1}(\cdot,x_i + hb^n_i) - V^n(\cdot,x_i) +\cfrac{1}{2}h(\sigma^n_i)^2(\mathcal{D}^2V^n)_i+hf^n_i )+ 
       &j ((H^{c}V^{n-1})_i -V^n_i+ \omicron_i^n ) \bigg] + \\
        &i   ((H^{\chi}V^{n-1})_i -V^n_i + \omicron_i^n  ) \Bigg\} =0,
\end{aligned}
\\ \\
 V^0_i=g(x_i).
\end{array}
\right .
\end{split}
\end{align}

\end{pro}
\begin{proof}
First we introduce the property below that will be useful for our demonstration. \\
 For any positive number $a,b,a^{'}$ and $b^{'}$ , the equation $\max \{ \min[A,B],C\} =0 $ is equivalent to the equation : 
   \begin{equation}\label{maxminlabel}
       \max \limits_{i\in \{0,1\}} \bigg\{(1-i)a \hspace{0.2cm} \min\limits_{j\in \{ 0,1\}} \big[ (1-j)a^{'} A + jb^{'} B\big]  + ibC  \bigg \} =0.
   \end{equation}
   The same conclusion holds for the inequalities :
   $$
  \max \{ \min[A,B],C\} \geq  0 \qquad  \text{and} \qquad \max \{ \min[A,B],C\} \leq 0.
   $$
Using (\ref{maxminlabel}), with $a =a^{'} = b =b^{'} = 1$, we rewrite (\ref{lagrange11st})  as the follow 
\begin{align*}
\begin{split}
\max\limits_{i \in \{0,1\}} 
       \Bigg\{ (1-i) \min \limits_{j\in \{ 0,1\}} \bigg[-(1-j)(\tilde{V}^{n-1}(\cdot,x_i + hb^n_i) - V^n(\cdot,x_i) +\cfrac{1}{2}h(\sigma^n_i)^2(\mathcal{D}^2V^n)_i+hf^n_i )+
       &j ( V^n_i - (H^{c}V^{n-1})_i  - \omicron_i^n  )  \bigg] + \\
        &i   (V^n_i - (H^{\chi}V^{n-1})_i   - \omicron_i^n  ) \Bigg\} =0,
    \end{split}
\end{align*}
thus, it follows 
\begin{align*}
\begin{split}
\max\limits_{i \in \{0,1\}} 
       \Bigg\{ -(1-i) \max \limits_{j\in \{ 0,1\}} \bigg[(1-j)(\tilde{V}^{n-1}(\cdot,x_i + hb^n_i) - V^n(\cdot,x_i) +\cfrac{1}{2}h(\sigma^n_i)^2(\mathcal{D}^2V^n)_i+hf^n_i )+
       &j (  (H^{c}V^{n-1})_i -V^n_i   + \omicron_i^n )  \bigg] + \\
        &i   (V^n_i - (H^{\chi}V^{n-1})_i   -\omicron_i^n   ) \Bigg\} =0,
    \end{split}
\end{align*}
from which we deduce 
\begin{align}\label{demo}
\begin{split}
\min\limits_{i \in \{0,1\}} 
       \Bigg\{ (1-i) \max \limits_{j\in \{ 0,1\}} \bigg[(1-j)(\tilde{V}^{n-1}(\cdot,x_i + hb^n_i) - V^n(\cdot,x_i) +\cfrac{1}{2}h(\sigma^n_i)^2(\mathcal{D}^2V^n)_i+hf^n_i )+
       &j (  (H^{c}V^{n-1})_i -V^n_i   +\omicron_i^n )  \bigg] + \\
        &i   ( (H^{\chi}V^{n-1})_i -V^n_i  +\omicron_i^n   ) \Bigg\} =0.
    \end{split}
\end{align}
\end{proof}
It is understood that  equation (\ref{demo}) holds for $n = 1, \dots, N$ and $i= 0, \dots ,M$.

Now our goal is to  show that the solution $V^h$ of our numerical scheme converges to the viscosity solution of the HJBQVI (\ref{eq:HJBI}). As introduced in the beginning, this  require to show that the numerical scheme $S$ is monotone, stable and non-locally consistent so that its unique solution  $V^h$ converges to the viscosity  solution $V$ of the HJBQVI as $h$ goes to  $0$. In order to prove the convergence, we recall the assumptions (\ref{ass1}) as well the proposition (\ref{propopo}) that will be useful for our analysis. Mainly the assumption $[\textbf{H}_{c,\chi}]$ which can be interpreted as the player paying a positive cost $(\underset{[0,T]\times \mathcal{U}}{inf}c \geq k,\underset{[0,T]\times \mathcal{V}}{inf}\chi\geq k)$ for having the right to perform his impulse.
\subsubsection{Monotonicity} 

\begin{lem} The  scheme $S$ is monotone. 
\end{lem}
\begin{proof}
We check the monotonicity property at grid points.  We assume $V$ and $\hat{V}$ are two functions such that $V \leq \hat{V}$ and  $V^n_i = \hat{V}^n_i$ at an arbitrary grid point $(\tau^n,x_i)$. Using the  definition, $S$ is monotone if $S(\cdot,\cdot,\cdot,V,\cdot,\cdot)  - S(\cdot,\cdot,\cdot,\hat{V},\cdot,\cdot) \geq 0$.\\
by (\ref{proposition31}), if $n=0$
$\qquad
S(\cdot,\cdot,\cdot,V,\cdot,\cdot)  - S(\cdot,\cdot,\cdot,\hat{V},\cdot,\cdot) = V^0_i - \hat{V}^0_i =  0.
$\\
If $n>0$,  we define $\ell:= (H^{c}V^{n-1})_i$ and $\ell^{'}:=( H^{\chi}V^{n-1})_i$ , by (\ref{demo})
\begin{align*}
 S(\cdot,\cdot&,\cdot,V,\cdot,\cdot)  - S(\cdot,\cdot,\cdot,\hat{V},\cdot,\cdot) \\
 \geq \min \limits_{i \in \{0,1\}} \Bigg\{ 
 & (1-i) \max\limits_{j\in \{ 0,1\}} \bigg[(1-j)(\tilde{\hat{V}}^{n-1}(\cdot,x_i + hb^n_i) - \hat{V}^n(\cdot,x_i) +\cfrac{1}{2}h(\sigma^n_i)^2(\mathcal{D}^2\hat{V}^n)_i+hf^n_i) + j  (\ell -\hat{V}^n_i+ \omicron_i^n  )  \bigg]  \\
        -& (1-i) \max\limits_{j\in \{ 0,1\}} \bigg[(1-j)(\tilde{V}^{n-1}(\cdot,x_i + hb^n_i) - V^n(\cdot,x_i) +\cfrac{1}{2}h(\sigma^n_i)^2(\mathcal{D}^2V^n)_i+hf^n_i) +  j  (\ell -V^n_i+ \omicron_i^n )  \bigg] \\
        & +i   (V^n_i -\hat{V}^n_i + \hat{\ell^{'}} -\ell^{'} + \omicron_i^n -\omicron_i^n   ) \Bigg \}, \\
\geq \min \limits_{i \in \{0,1\}} \Bigg\{ 
 & (1-i) \min\limits_{j\in \{ 0,1\}} \bigg[ (1-j)(\tilde{\hat{V}}^{n-1}(\cdot,x_i + hb^n_i ) - \hat{V}^n(\cdot,x_i)) +\cfrac{1}{2}h(\sigma^n_i)^2(\mathcal{D}^2\hat{V}^n)_i+hf^n_i)   \\
        & \qquad \qquad \qquad   -(1-j)(\tilde{V}^{n-1}(\cdot,x_i + hb^n_i) - V^n(\cdot,x_i) +\cfrac{1}{2}h(\sigma^n_i)^2(\mathcal{D}^2V^n)_i+hf^n_i) \\
        &  +j   ( \omicron_i^n   -\omicron_i^n  )\bigg] +i   ( \omicron_i^n   -\omicron_i^n    ) \Bigg \}, \\
  \geq \min \limits_{i \in \{0,1\}} \Bigg\{ 
  & (1-i) \min\limits_{j\in \{ 0,1\}} \bigg[ (1-j)(\tilde{\hat{V}}^{n-1}(\cdot,x_i + hb^n_i) -\tilde{V}^{n-1}(\cdot,x_i + hb^n_i)  + \cfrac{1}{2}h\sigma(\cdot)^2[\mathcal{D}^2U^n]_i ) 
  \bigg]  
 \Bigg\}  ,
\end{align*}
from what we have 
\begin{align}\label{demono}
\begin{split}
 S(\cdot,\cdot,\cdot,V,\cdot,\cdot)  - S(\cdot,\cdot,\cdot,\hat{V},\cdot,\cdot)  
    \geq \min \limits_{i \in \{0,1\}} \Bigg\{ 
   &(1-i) \min\limits_{j\in \{ 0,1\}} \bigg[ (1-j)(\tilde{\hat{V}}^{n-1}(\cdot,x_i + hb^n_i) -\tilde{V}^{n-1}(\cdot,x_i + hb^n_i))  + \cfrac{1}{2}h\sigma(\cdot)^2[\mathcal{D}^2U^n]_i 
  \bigg], 
  \end{split}
\end{align}
and 
\begin{align*}
  \tilde{\hat{V}}^{n-1}(\cdot,x_i + b^n_i\Delta t) &
  -\tilde{V}^{n-1}(\cdot,x_i + b^n_i\Delta t)   \\&:=  \alpha \Big( \hat{V}(\tau^n,x_{k+1}) - V(\tau^n,x_{k+1}) \Big)  +  (1-\alpha)\Big( \hat{V}(\tau^n,x_k)-V(\tau^n,x_k)  \Big ) \geq 0. 
\end{align*}
We defined the vector $U$ as $U := \hat{V}^n - V^n  = ( \hat{V}^n_0 -V^n_0 , \cdots, \hat{V}^n_M -V^n_M )$ with positive entries.\\
 Then we have  $\mathcal{D}^2U :=\dfrac{U_{i+1}  + U_{i-1}  }{(\Delta x)^2}\geq 0$, 
 and the right hand component of (\ref{demono}) is positive hence $S(\cdot,\cdot,\cdot,V,\cdot,\cdot)  - S(\cdot,\cdot,\cdot,\hat{V},\cdot,\cdot) \geq 0$ as desired.
Hence,  we obtain $S(\cdot,\cdot,\cdot,V,\cdot,\cdot)  - S(\cdot,\cdot,\cdot,\hat{V},\cdot,\cdot) \geq 0$ as desired.
\end{proof}
\subsubsection{Stability}
\begin{lem}The  scheme $S$ is stable. 
\end{lem}
\begin{proof}
Let $V = V^h$ be a solution  of the scheme $S$. We have by (\ref{proposition31}) if $n =0$,
$$
 0  = S\Big(h,(t,x),V^h, \cdot , \cdot \Big) = V_i^0 - g(x_i),
$$
   which implies  
\begin{equation} \label{boundvo}
-||g||_{\infty} \leq  V_i^0 \leq  ||g||_{\infty} . 
\end{equation}
If $n > 0$, by (\ref{demo}) we have 
\begin{align*}
    \begin{split}
       0 &= S\Big(h,(t,x),V^h,  \cdot , \cdot  \Big)  \\
 & \leq \max \limits_{j\in \{ 0,1\}} \bigg[(1-j)(\tilde{V}^{n-1}(\cdot,x_i + hb^n_i) - V^n(\cdot,x_i) +\cfrac{1}{2}h(\sigma^n_i)^2(\mathcal{D}^2V^n)_i+hf^n_i )+
       j (  (H^{c}V^{n-1})_i -V^n_i   + \cfrac{1}{2}h(\sigma^n_i)^2(\mathcal{D}^2V^n)_i )  \bigg] \\
 &  = \max \limits_{j\in \{ 0,1\}} \bigg[(1-j)(\tilde{V}^{n-1}(\cdot,x_i + hb^n_i) - V^n(\cdot,x_i) +hf^n_i )+
       j (  (H^{c}V^{n-1})_i -V^n_i )  + \cfrac{1}{2}h(\sigma^n_i)^2(\mathcal{D}^2V^n)_i   \bigg]\\
 &  = \max  \bigg\{\tilde{V}^{n-1}(\cdot,x_i + hb^n_i) - V^n(\cdot,x_i) +hf^n_i ,
         (H^{c}V^{n-1})_i -V^n_i  \bigg\} + \cfrac{1}{2}h(\sigma^n_i)^2(\mathcal{D}^2V^n)_i  ,
    \end{split}
\end{align*}
where we have employed the fact that  $ \omicron_i^n := \cfrac{1}{2}h(\sigma^n_i)^2(\mathcal{D}^2V^n)_i $ (see remark \ref{remarkomicron}).\\
Now let's pick $j$ such that $V^n_j = \max \limits_{k} V^n_k$ from which it  follows  that $V^n_j \geq V^n_{j \pm 1}$ such that $[\mathcal{D}^2V^n]_j \leq 0$ and 
\begin{equation}\label{setij}
0 \leq  \max  \bigg\{\tilde{V}^{n-1}(\cdot,x_i + hb^n_i) - V^n(\cdot,x_i) +hf^n_i ,
         (H^{c}V^{n-1})_i -V^n_i  \bigg\} .
\end{equation}
Therefore, setting $i=j$ in (\ref{setij}) and taking each component separately  we have 
\begin{equation}\label{stab1c}
\tilde{V}^{n-1}(\cdot,x_i + hb^n_i) - V^n(\cdot,x_i) +hf^n_i \leq  \max \limits_{k} V^{n-1}_k - V^n_i + h ||f||_{\infty},
\end{equation}
and 
\begin{equation}\label{stab2c}
(H^{c}V^{n-1})_i = \sup \limits_{ (x+\xi) \in \mathcal{U}^h}[\tilde{V}(\tau^n, x_i + \xi)-c(\tau^n, \xi)] < \sup \limits_{ (x+\xi) \in  \mathcal{U}^h} [\tilde{V}(\tau^n, x_i + \xi)] \leq 
 \max \limits_{k} V^{n-1}_k.
\end{equation}
From which we deduce that 
\begin{align*}
\begin{split}
0 \leq \max \limits_{j\in \{ 0,1\}} \bigg[(1-j)(\max \limits_{k} V^{n-1}_k - V^n_i + h ||f||_{\infty})+
       j &(  \max \limits_{k} V^{n-1}_k -V^n_i )     \bigg]\\
       & = \max \limits_{k} V^{n-1}_k - V^n_i + h ||f||_{\infty}),
       \end{split}
\end{align*}
hence 
\begin{align}\label{equiva}
 V^n_i \leq \max \limits_{k} V^{n-1}_k  + h ||f||_{\infty}.
\end{align}
Combining (\ref{equiva}) and (\ref{boundvo}), we obtain 
\begin{align}\label{equiva2}
\max \limits_{k} V^{n}_k  \leq   nh ||f||_{\infty} + ||g||_{\infty} \leq  ||f||_{\infty}T + ||g||_{\infty} \qquad \text{for all } n,
\end{align}
which provide an upper bound for  the solution $V^h$. \\
For the lower bound, let $n >0$ and pick $j$ such that $V^n_j = \min \limits_{k} V^n_k$  from which it  follows  that $V^n_j \leq V^n_{j \pm 1}$ such that $[\mathcal{D}^2V^n]_j \geq 0$.
\begin{align*}
\begin{split}
0 &= S\Big(h,(t,x),V^h,  \cdot , \cdot  \Big) \\
 &= \min\limits_{i \in \{0,1\}} 
       \Bigg\{ (1-i) \max \limits_{j\in \{ 0,1\}} \bigg[(1-j)(\tilde{V}^{n-1}(\cdot,x_i + hb^n_i) - V^n(\cdot,x_i) +\cfrac{1}{2}h(\sigma^n_i)^2(\mathcal{D}^2V^n)_i+hf^n_i )+ 
       j (  (H^{c}V^{n-1})_i -V^n_i   + \\ 
        &  \qquad \qquad \qquad \qquad \qquad \cfrac{1}{2}h(\sigma^n_i)^2(\mathcal{D}^2V^n)_i )  \bigg]  + i ( (H^{\chi}V^{n-1})_i -V^n_i  +\cfrac{1}{2}h(\sigma^n_i)^2(\mathcal{D}^2V^n)_i  ) \Bigg\} ,
    \end{split}
\end{align*}
hence , at least one of 
\begin{align}\label{oneof1}
\begin{split}
0 = \max \limits_{j\in \{ 0,1\}} \bigg[(1-j)(\tilde{V}^{n-1}(\cdot,x_i + hb^n_i) - V^n(\cdot,x_i) +\cfrac{1}{2}h(\sigma^n_i)^2(\mathcal{D}^2V^n)_i+hf^n_i )+ 
       j (  (H^{c}V^{n-1})_i -V^n_i   + \cfrac{1}{2}h(\sigma^n_i)^2(\mathcal{D}^2V^n)_i )  \bigg],
    \end{split}
\end{align}
or 
\begin{align}\label{oneof2}
\begin{split}
0 = (H^{\chi}V^{n-1})_i -V^n_i  +\cfrac{1}{2}h(\sigma^n_i)^2(\mathcal{D}^2V^n)_i   ,
    \end{split}
\end{align}

holds necessarily. 
From (\ref{oneof2}) we have 
\begin{align*}
    \begin{split}
        (H^{\chi}V^{n-1})_i -V^n_i &= \inf \limits_{(x+
 \eta)\in \mathcal{V^h}}[\tilde{V}(\tau^n,x+\eta) +\chi(\tau^n,\eta)].\\
 & > \inf \limits_{(x+
 \eta)\in \mathcal{V^h}}[\tilde{V}(\tau^n,x+\eta)] \geq   \min \limits_{k} V^{n-1}_k,
    \end{split}
\end{align*}
 and 
 \begin{align*}
    \begin{split}
        0  \geq   \min \limits_{k} V^{n-1}_k -V^n_i  +\cfrac{1}{2}h(\sigma^n_i)^2(\mathcal{D}^2V^n)_i \geq  \min \limits_{k} V^{n-1}_k  -V^n_i .
    \end{split}
\end{align*}

From (\ref{oneof1}), we have  
\begin{align*}
\begin{split}
0&= \max  \bigg\{(\tilde{V}^{n-1}(\cdot,x_i + hb^n_i) - V^n(\cdot,x_i) +\cfrac{1}{2}h(\sigma^n_i)^2(\mathcal{D}^2V^n)_i+hf^n_i ),
          (H^{c}V^{n-1})_i -V^n_i   + \cfrac{1}{2}h(\sigma^n_i)^2(\mathcal{D}^2V^n)_i )  \bigg\} \\
&\geq (\tilde{V}^{n-1}(\cdot,x_i + hb^n_i) - V^n(\cdot,x_i) +\cfrac{1}{2}h(\sigma^n_i)^2(\mathcal{D}^2V^n)_i+hf^n_i )  \\
& \geq  (\tilde{V}^{n-1}(\cdot,x_i + hb^n_i) - V^n(\cdot,x_i) +hf^n_i \geq   \min \limits_{k} V^{n-1}_k  -V^n_i - h||f||_{\infty}.
    \end{split}
\end{align*}

From which we  deduce that 
\begin{align}\label{deduce}
    \begin{split}
     V^n_i \geq   \min \limits_{k} V^{n-1}_k  - h||f||_{\infty}.
    \end{split}
\end{align}
Combining (\ref{deduce}) with (\ref{boundvo})  we obtain 
\begin{align}\label{equiva3}
\min \limits_{k} V^{n}_k  \geq  - nh ||f||_{\infty} - ||g||_{\infty} \geq  -||f||_{\infty}T  -||g||_{\infty} \qquad \text{for all } n,
\end{align}
providing a lower bound for the solution $V$.
\end{proof}

\subsubsection{Non local consistency}
To prove the consistency of our scheme, we need before to ensure the convergence of the finite non-empty sets $\mathcal{U}^h$ and $\mathcal{V}^h$ to the original set which they approximate. This can be obtain using the Hausdorff metric.   We establish the proposition below for this matter. 
\begin{pro} \label{Haudsorff}
    The approximated sets $\mathcal{U}^h$ (resp., $\mathcal{V}^h$) converges locally uniformly (with respect to the Hausdorff metric) to $\mathcal{U}$ (resp., $\mathcal{V}$) as $h\to 0$. Moreover $(t,x) \mapsto \mathcal{U}^h (t,x)	$ and $(t,x) \mapsto \mathcal{V}^h (t,x)	$ are continuous (with respect to the Hausdorff metric) for each $h$.
\end{pro}
Below, we state few lemmas that are also useful for our consistency proof. 
\begin{lem} \label{consistency1}
Let $(a_m)_m$, $(b_m)_m$, $(c_m)_m$ be sequences of real numbers such that  for each $ m$ ,  $c_m \geq \min\{a_m,b_m\}$ (resp., $c_m \leq \min\{a_m,b_m\}$). \\
Then, 
\begin{align} \label{consitencylabel1}
    &\liminf \limits_{ m \to \infty} c_m \geq \min\{ \liminf \limits_{ m \to \infty} a_m, \liminf \limits_{ m \to \infty} b_m \}  ,
    \\
    (\text{resp.  } &\limsup \limits_{ m \to \infty} c_m \leq \min\{ \limsup \limits_{ m \to \infty} a_m, \limsup \limits_{ m \to \infty} b_m \}).
\end{align}
Similarly, let $(a_m)_m$, $(b_m)_m$, $(c_m)_m$ be sequences of real numbers such that  for each $ m$ ,  $c_m \geq \max\{a_m,b_m\}$ (resp., $c_m \leq \max\{a_m,b_m\}$). \\
Then, 
\begin{align} \label{consitencylabel2}
    &\liminf \limits_{ m \to \infty} c_m \geq \max\{ \liminf \limits_{ m \to \infty} a_m, \liminf \limits_{ m \to \infty} b_m \}  ,
\end{align}
\begin{align}\label{consistencylabel22}
    (\text{resp.  } &\limsup \limits_{ m \to \infty} c_m \leq \max\{ \limsup \limits_{ m \to \infty} a_m, \limsup \limits_{ m \to \infty} b_m \}).
\end{align}
\end{lem}
\begin{proof} We start by proving (\ref{consitencylabel1}).
    The proof follows the assumption that  for each $m$, if $c_m \geq \min\{a_m,b_m\}$ then
    $$
    \inf \limits_{k\geq m} c_k \geq \inf \limits_{k\geq m} \min\{ a_k,b_k \} = \min \{ 
 \inf \limits_{k\geq m} a_k, \inf \limits_{k\geq m} b_k\},
    $$
    therefore, taking limits into the above and using the continuity of the function $(x,y) \mapsto \min(x,y)$ coupled with  $\lim \limits_{ m \to \infty} \inf \limits_{k\geq m} x_k =  \liminf \limits_{ m \to \infty}  x_m$ we obtain 
    $$
    \lim \limits_{ m \to \infty} \inf \limits_{k\geq m} c_k \geq  \lim \limits_{ m \to \infty} \min \{ 
 \inf \limits_{k\geq m} a_k, \inf \limits_{k\geq m} b_k\} =  \min \{ 
 \lim \limits_{ m \to \infty} \inf \limits_{k\geq m} a_k, \lim \limits_{ m \to \infty} \inf \limits_{k\geq m} b_k\}
     = \min \{ 
 \liminf\limits_{ m \to \infty}  a_m, \liminf \limits_{ m \to \infty}  b_m\}.$$
 Using analogous argument we obtain the results for  $c_m \leq \min\{a_m,b_m\}$ where 
$$
    \sup \limits_{k\geq m} c_k \leq \sup \limits_{k\geq m} \min\{ a_k,b_k \} \leq \min \{ 
 \sup \limits_{k\geq m} a_k, \sup \limits_{k\geq m} b_k\}.
    $$
 Similarly for (\ref{consitencylabel2}), if  $c_m \geq \max\{a_m,b_m\}$, then 
 $$
 \inf \limits_{k\geq m} c_k \geq \inf \limits_{k\geq m} \max\{ a_k,b_k \} = \max \{ 
 \inf \limits_{k\geq m} a_k, \inf \limits_{k\geq m} b_k\},
 $$
  therefore, taking limits into the above and using the continuity of the function $(x,y) \mapsto \max(x,y)$ coupled with  $\lim \limits_{ m \to \infty} \inf \limits_{k\geq m} x_k =  \liminf \limits_{ m \to \infty}  x_m$ we obtain 
  $$
    \lim \limits_{ m \to \infty} \inf \limits_{k\geq m} c_k \geq  \lim \limits_{ m \to \infty} \max \{ 
 \inf \limits_{k\geq m} a_k, \inf \limits_{k\geq m} b_k\} =  \max \{ 
 \lim \limits_{ m \to \infty} \inf \limits_{k\geq m} a_k, \lim \limits_{ m \to \infty} \inf \limits_{k\geq m} b_k\}
     = \max \{ 
 \liminf\limits_{ m \to \infty}  a_m, \liminf \limits_{ m \to \infty}  b_m\}.$$
 Using analogous argument we obtain the results for  $c_m \leq \max\{a_m,b_m\}$ where 
 $$
    \sup \limits_{k\geq m} c_k \leq \sup \limits_{k\geq m} \min\{ a_k,b_k \} \leq \min \{ 
 \sup \limits_{k\geq m} a_k, \sup \limits_{k\geq m} b_k\}.
    $$
\end{proof}
Next, we must ensure the approximated operators are suited for the intervention operators.
\begin{lem} \label{consistency2}
Let $(v^h)_{h > 0}$ be a family of uniformly bounded real-valued maps from $[O,T]\times \mathbb{R}$ wit half-relaxed limits $\overline{v}$ and $\underline{v}$. Let  $v^{n,h} = ( v^h(\tau^n, x_0), \dots, v^h(\tau^n,x_M))^{\intercal}$ be a vector  whose components are obtained by evaluating the function $v^h(\tau^n,\cdot)$ for each $x$ on the numerical grid.  Then, for any $(t,x) \in [0,T]\times \mathbb{R}$,
\begin{equation} \label{consistency2.1}
  H_{sup}^c \underline{v}(t,x) \leq \liminf\limits_{ \substack{h\to 0 \\ (\tau^n, x_i) \to (t,x)}} (H^c v^{n,h})_i \leq  \limsup\limits_{ \substack{h\to 0 \\ (\tau^n, x_i) \to (t,x)}} (H^c v^{n,h})_i \leq H_{sup}^c \overline{v}(t,x)  ,
\end{equation}

and 
\begin{equation} \label{consistency2.2}
  H_{inf}^{\chi} \underline{v}(t,x) \leq \liminf\limits_{ \substack{h\to 0 \\ (\tau^n, x_i) \to (t,x)}} (H^{\chi} v^{n,h})_i \leq  \limsup\limits_{ \substack{h\to 0 \\ (\tau^n, x_i) \to (t,x)}} (H^{\chi} v^{n,h})_i \leq H_{inf}^{\chi} \overline{v}(t,x)  .
\end{equation}

\end{lem}
\begin{proof} We start by proving the leftmost inequality of (\ref{consistency2.1}). 
    Let $(h_m, s_m,y_m)$ be sequence that satisfies  $h_m \to 0$ and $(s_m,y_m) \to (t,x)$ as $m\to 0$. The choice of $s_m$ and $y_m$ are made wisely so that $s_m  = \tau^{n_m} = T- n_m\Delta t$ and $y_m = x_{im}$. Let $\delta > 0$ and $\xi^{\delta} \in \mathcal{U}(t,x)$ such that 
    \begin{equation}\label{440}
         H_{sup}^c \underline{v}(t,x) = \sup \limits_{\xi \in \mathcal{U}}[ \underline{v}(t, x + \xi) - c(t,\xi)] \leq \underline{v}(t,x + \xi^{\delta})- c(t,\xi^{\delta}) + \delta.
    \end{equation}

    Using  Proposition (\ref{Haudsorff}), we can pick a sequence $(\xi_m)_m$ such that 
    $xi_m \to \xi^{\delta}$ and  $xi_m \in \mathcal{U}^{h_m}(t,x)$ for each $m$. We have 
    \begin{align*}
    (H^{c} v^{n_m,h_m})_{i_m} &= \sup \limits_{\xi \in \mathcal{U}^{h_m}}\big[\tilde{v}^{h_m}(s_m,y_m + \xi ) - c(s_m, \xi)\big]\\
    &\geq \tilde{v}^{h_m}(s_m,y_m + \xi_m ) - c(s_m, \xi_m) \\
    & = \alpha_m v^{h_m} (s_m,x_{k_m +1}) + (1- \alpha_m ) v^{h_m} (s_m,x_{k_m}) - c(s_m, \xi_m) \qquad 0\leq \alpha_m \leq 1,
    \end{align*} 
    with $k_m$ satisfying $x_{k_m} \to x+ \xi_{\delta}$ as $m\to \infty$ and $(s_m, x_{k_m}) \to (t, x+ \xi_{\delta} )$ as $m \to \infty$.
    Using Lemma (\ref{consistency1}), 
    \begin{align}\label{443}
    \begin{split}
       \liminf \limits_{ m\to \infty} (H^{c} v^{n_m,h_m})_{i_m}  &\geq  \liminf \limits_{ m\to \infty}\{ \min \{v^{h_m} (s_m,x_{k_m +1}), v^{h_m} (s_m,x_{k_m})   \} - c(s_m,\xi_m)\} \\
       & \geq \min \{ \liminf \limits_{ m\to \infty} v^{h_m} (s_m,x_{k_m +1}), \liminf \limits_{ m\to \infty} v^{h_m} (s_m,x_{k_m}) \} - c(s_m,\xi_m).
       \end{split}
    \end{align}
Therefore, by the definition of the half relaxed limit $\underline{v}$, 
\begin{equation} \label{444}
    \liminf \limits_{ m\to \infty} v^{h_m} (s_m,x_{k_m}) \geq \liminf \limits_{ \substack{h\to 0 \\ (s,y) \to (t, x+ \xi_{\delta} )}}  v^h(s,y) = \underline{v}(t, x+ \xi_{\delta} ) ,
\end{equation}
and similarly, since $x_{k_m +1} - x_{k_m} = \Delta x$, 
\begin{equation}\label{445}
    \liminf \limits_{ m\to \infty} v^{h_m} (s_m,x_{k_m+1}) \geq  \underline{v}(t, x+ \xi_{\delta} ) .
\end{equation}
Applying (\ref{440}), (\ref{444}) and (\ref{445}) to (\ref{443}), we obtain 
\begin{equation} \label{446}
    \liminf \limits_{ m\to \infty} (H^{c} v^{n_m,h_m})_{i_m} \geq \underline{v}(t, x+ \xi_{\delta} )  - c(t, x + \xi_{\delta}) \geq  H_{sup}^c \underline{v}(t,x) -\delta,
\end{equation}
and since $\delta$ is arbitrary, we obtain the leftmost inequality. \\
Now, we prove the rightmost inequality of (\ref{consistency2.1}). We use the same sequence $(h_m, s_m, y_m)$ as in the previous paragraph. Since the defined set $\mathcal{U}^{h_m}(s_m,y_m)$ is finite and compact by induction for each $m$, there exists a $\xi_m \in \mathcal{U}^{h_m}(s_m,y_m) $  such that 
\begin{align*}
   (H^{c} v^{n_m,h_m})_{i_m}   &= \sup\limits_{\xi \in \mathcal{U}^{h_m}(s_m,y_m)  } \big[
    \tilde{v}^{h_m}(s_m,y_m + \xi ) - c(s_m, \xi) 
   \big] \\
   & =  \tilde{v}^{h_m}(s_m,y_m + \xi_m ) - c(s_m, \xi_m) \\
   & = \alpha_m v^{h_m} (s_m,x_{k_m +1}) (1- \alpha_m ) v^{h_m} (s_m,x_{k_m}) - c(s_m, \xi_m),
\end{align*}
where $ 0 \leq \alpha_m \leq 1$ with  $k_m$ satisfying $x_{k_m} \to x+ \xi_{\delta}$ as $m\to \infty$.  Using proposition (\ref{Haudsorff}), the sequence $(\xi_m)_m$ is contained in a compact set and converges to a subsequence $(\xi_{m_j})_j$  with limit $\hat{\xi}$. Hence the right hand side of the inequality 
$$
d(\hat{\xi}, \mathcal{U}(t,x)) \leq d(\hat{\xi},\xi_{m_j}) + d(\xi_{m_j},\mathcal{U}(t,x)) ,
$$
converges to 0 as $j\to 0$ and $\hat{\xi} \in \mathcal{U}(t,x)$. \\
Similarly to (\ref{443}) and (\ref{446}), 

\begin{align*} 
\begin{split}
       \limsup \limits_{ m\to \infty} (H^{c} v^{n_m,h_m})_{i_m}  &\leq  \limsup\limits_{ m\to \infty}\{ \max \{v^{h_m} (s_m,x_{k_m +1}), v^{h_m} (s_m,x_{k_m})   \} - c(s_m,\xi_m)\} \\
       & \leq \max \{ \limsup \limits_{ m\to \infty} v^{h_m} (s_m,x_{k_m +1}), \limsup\limits_{ m\to \infty} v^{h_m} (s_m,x_{k_m}) \} - c(s_m,\xi_m)
       \end{split} \\ 
       & \leq \overline{v}(t, x+ \hat{\xi})  - c(t,x+\hat{\xi}) \\ 
       & \leq H_{sup}^c \overline{v}(t,x).
    \end{align*}
    And we obtain the result  for (\ref{consistency2.1}). 
    We now prove (\ref{consistency2.2}) with the same reasoning.  
\end{proof}

We are now ready to establish the nonlocal consistency of our numerical scheme.

\begin{lem}  The  scheme $S$ is non-locally consistent.
\end{lem}
\begin{proof}
Let $\Omega = [0,T) \times \mathbb{R}$ and $\varphi \in C^{1,2}(\overline{\Omega})$ such that $\varphi^n = ( \varphi(\tau^n, x_0), \dots, \varphi(\tau^n,x_M))^{\intercal}$. For simplicity, we  let $\varphi^n_i = \varphi(\tau^n, x_i)$.  Let $(V^h)_{h>0}$ be a family of uniformly bounded real-valued maps from $\overline{\Omega}$ with half-relaxed limits $\overline{V}$ and $\underline{V}$.
 Let $(h_m, s_m,y_m,\xi_m)$ be an arbitrary sequence that satisfies  $h_m \to 0$ and $(s_m,y_m)\to (t,x)$ and $\xi_m \to 0$ as $m\to 0$. We assume the sequence is chosen such that $s_m  = \tau^{n_m} = T- n_m\Delta t$ and $y_m = x_{im}$. By the definition (\ref{lagrange1st}) of the scheme $S$, 
 \begin{align} \label{448}
\begin{split}
S \Big(h_m,(s_m,y_m),\varphi + \xi_m , \mathcal{H}^{h_m}_{sup}V^{h_m},\mathcal{H}^{h_m}_{inf} V^{h_m}\Big) = \left \{
\begin{array}{l}
 \max\big \{\min[S^{(1)}_m,S^{(2)}_m], S^{(3)}_m\big \} \qquad  \text{ if } n_m >0,\\ 
 S^{(4)}_m \qquad \qquad \qquad \qquad \qquad \quad\text{        if } n_m  =0.
\end{array}
\right .
\end{split}
\end{align} 
where 
\begin{align} \label{449}
\begin{split}
S^{(1)}_m &= -\dfrac{(\tilde{(\varphi + \xi_m)}^{n_m-1}(\cdot,y_m + h_mb^{n_m}_{i_m}) - \varphi^{n_m}(\cdot,y_m) -\xi_m)}{h_m} -\cfrac{1}{2}(\sigma^{n_m}_{i_m})^2(\mathcal{D}^2_{xx}\varphi^{n_m})_{i_m}-f^{n_m}_{i_m}   , \\
S^{(2)}_m &=  \varphi^{n_m}_{i_m} + \xi_m-(H^{c} V^{n_m-1,h_m})_{i_m} ,
\\
S^{(3)}_m &=  \varphi^{n_m}_{i_m} + \xi_m-(H^{\chi}V^{n_m-1,h_m})_{i_m},
\\
S^{(4)}_m &= \varphi^{n_m}_{i_m} + \xi_m  - g(y_m).
\end{split}
\end{align} 
Firstly, we have by a Taylor expansion 
\begin{align}\label{461}
\begin{split}
&\dfrac{(\tilde{(\varphi + \xi_m)}^{n_m-1}(\cdot,y_m + \Delta_{\tau} b^{n_m}_{i_m}) - \varphi^{n_m}(\cdot,y_m) -\xi_m)}{\Delta_{\tau}}=\dfrac{\varphi(s_m + \Delta_{\tau}, y_m +  \Delta_{\tau} b^{n_m}_{i_m}) - \varphi^{n_m}(\cdot,y_m) + \xi_m-\xi_m }{\Delta_{\tau}} + O(\dfrac{(\Delta_{\tau})^2}{\Delta_x}) \\
&\qquad = \varphi_t(t,x) + b(t,x)\varphi_x(t,x) + O(\dfrac{(\Delta_{\tau})^2}{\Delta_x} + \Delta_{\tau}) = \varphi_t(t,x) + b(t,x)\varphi_x(t,x) + O(h_m),
\end{split}
\end{align} 
where $\varphi_t$ (resp, $\varphi_x$) is the first order time derivative (resp, first order space derivative) of the function $\varphi$. Then
\begin{align}\label{450}
    \begin{split}
\lim \limits_{m \to \infty} S^{(1)}_m &= -\lim \limits_{m \to \infty} \Big \{ \dfrac{(\tilde{(\varphi + \xi_m)}^{n_m-1}(\cdot,y_m + h_mb^{n_m}_{i_m}) - \varphi^{n_m}(\cdot,y_m) -\xi_m)}{h_m} +\cfrac{1}{2}(\sigma^{n_m}_{i_m})^2(\mathcal{D}^2_{xx}\varphi^{n_m})_{i_m}+f^{n_m}_{i_m} \Big \} \\ 
& = -\lim \limits_{m \to \infty} \Big \{ \varphi_t(t,x) + b^{n_m}_{i_m}\varphi_x(t,x) +  \cfrac{1}{2}(\sigma^{n_m}_{i_m})^2 \varphi_{xx}(t,x) +f^{n_m}_{i_m} +  O(h_m)   \Big \} \\
& = -\varphi_t(t,x)  - b(t,x)\varphi_x(t,x) -  \cfrac{1}{2}\sigma(t,x)^2 \varphi_{xx}(t,x)  - f(t,x)  \qquad \text{ as }\quad  O(h_m) \to 0 , m \to \infty.
    \end{split}
\end{align}
Now by using lemma (\ref{consistency2}), 
\begin{align}\label{451}
\begin{split}
    \liminf \limits_{m \to \infty} S^{(2)}_m &= \liminf \Big \{\varphi^{n_m}_{i_m} + \xi_m-(H^{c} V^{n_m-1,h_m})_{i_m} \Big \} \\
    &\geq   \varphi(t,x) - \limsup\limits_{  m \to \infty} (H^{c} V^{n_m-1,h_m})_{i_m} \geq \varphi(t,x)  - H_{sup}^c \overline{v}(t,x)  ,
    \end{split}
\end{align}
and 
\begin{align}\label{452}
\begin{split}
    \limsup \limits_{m \to \infty} S^{(2)}_m &= \limsup \Big \{\varphi^{n_m}_{i_m} + \xi_m-(H^{c} V^{n_m-1,h_m})_{i_m} \Big \} \\
    &\leq \varphi(t,x) -  \liminf \limits_{m \to \infty} (H^{c} V^{n_m-1,h_m})_{i_m} \leq \varphi(t,x)  - H_{sup}^c \underline{v}(t,x)  .
    \end{split}
\end{align}
Similarly 
\begin{align}\label{4511}
\begin{split}
    \liminf \limits_{m \to \infty} S^{(3)}_m &= \liminf \Big \{\varphi^{n_m}_{i_m} + \xi_m-(H^{\chi}V^{n_m-1,h_m})_{i_m} \Big \} \\
    &\geq   \varphi(t,x) - \limsup\limits_{  m \to \infty} (H^{\chi}V^{n_m-1,h_m})_{i_m} \geq \varphi(t,x)  - H_{inf}^{\chi} \overline{v}(t,x)  ,
    \end{split}
\end{align}
and 
\begin{align}\label{4522}
\begin{split}
    \limsup \limits_{m \to \infty} S^{(3)}_m &= \limsup \Big \{\varphi^{n_m}_{i_m} + \xi_m-(H^{\chi}V^{n_m-1,h_m})_{i_m} \Big \} \\
    &\leq \varphi(t,x) -  \liminf \limits_{m \to \infty} (H^{\chi}V^{n_m-1,h_m})_{i_m} \leq \varphi(t,x)  - H_{inf}^{\chi} \underline{v}(t,x)  .
    \end{split}
\end{align}

Now let's assume $t < T$, we may assume $s_m < T$ since $s_m \to T$. Hence, using Lemma (\ref{consistency1},eq. \ref{consitencylabel2})
\begin{align*}
    \begin{split}
        \liminf \limits_{m\to \infty}  S \Big(h_m,(s_m,y_m),\varphi + \xi_m , \mathcal{H}^{h_m}_{sup}V^{h_m},\mathcal{H}^{h_m}_{inf} V^{h_m}\Big) &=   \liminf \limits_{m\to \infty}  \max\big \{\min[S^{(1)}_m,S^{(2)}_m], S^{(3)}_m\big \} \\ 
        &\geq \max \big \{ \liminf \limits_{m\to \infty} \min[S^{(1)}_m,S^{(2)}_m], \liminf \limits_{m\to \infty} S^{(3)}_m  \big \} \\
        & \geq \max \big \{  \min \{ \liminf\limits_{m\to \infty} S^{(1)}_m, \liminf\limits_{m\to \infty}S^{(2)}_m\}, \liminf \limits_{m\to \infty} S^{(3)}_m  \big \}.
    \end{split}
\end{align*}
Applying (\ref{450}), (\ref{451}) and (\ref{4511})

\begin{align}\label{453}
    \begin{split}
        \liminf \limits_{m\to \infty}  S \Big(&h_m,(s_m,y_m),\varphi + \xi_m , \mathcal{H}^{h_m}_{sup}V^{h_m},\mathcal{H}^{h_m}_{inf} V^{h_m}\Big)\\ 
        &\geq  \max \Big\{ \min\big\{ -\varphi_t(t,x)  - b(t,x)\varphi_x(t,x) -  \cfrac{1}{2}\sigma(t,x)^2 \varphi_{xx}(t,x)  - f(t,x) , \varphi(t,x)  - H_{sup}^c \overline{v}(t,x)  \big \}, \\
        &  \qquad  \varphi(t,x)  - H_{inf}^{\chi} \overline{v}(t,x) \Big\} = F_{*}(t,x,\phi(t,x),D\phi(t,x),D^2\phi(t,x),
H^c_{sup} \overline{v}(t,x),H^\chi_{inf} \overline{v}(t,x)).
    \end{split}
\end{align}
Now since $(h_m, s_m,y_m,\xi_m)$ is an arbitrary sequence, (\ref{453}) implies 
\begin{equation}\label{454}
 \liminf \limits_{\substack{h\to 0 \\s,y\to t, x \\ \epsilon \to 0  }} S(h,(s,y),\phi +\epsilon,\mathcal{H}^{h}_{sup} V^h,\mathcal{H}^{h}_{inf} V^h) \geq  F_{*}(t,x,\phi(t,x),D\phi(t,x),D^2\phi(t,x),
H^c_{sup} \overline{v}(t,x),H^\chi_{inf} \overline{v}(t,x)),  
\end{equation}
which is exactly the nonlocal consistency inequality (\ref{436}).\\
Symmetrically using (\ref{consistencylabel22})  we can establish the nonlocal consistency inequality (\ref{437}). 
\begin{equation} \label{455}
 \limsup \limits_{\substack{h\to 0 \\s,y\to t, x \\ \epsilon \to 0  }} S(h,(s,y),\phi +\epsilon,\mathcal{H}^{h}_{sup} V^h,\mathcal{H}^{h}_{inf} V^h) \leq  F^{*}(t,x,\phi(t,x),D\phi(t,x),D^2\phi(t,x),
H^c_{sup} \underline{v}(t,x),H^\chi_{inf} \underline{v}(t,x)).
\end{equation}
Now, when $t = T$, since $H^c_{sup}g \leq g \leq H^\chi_{inf}g$, it follows that $g(y_m) = \min \big \{ \max \big \{ g(y_m),H^c_{sup}g(y_m) \big \},H^\chi_{inf}g(y_m) \big\}$  for each $m$. Therefore, 
\begin{align*}
    S^{(4)}_m &= \varphi^{n_m}_{i_m} + \xi_m  - \min \big \{ \max \big \{ g(y_m),H^c_{sup}g(y_m) \big \},H^\chi_{inf}g(y_m) \big\}\\
    & = \max \big \{ \min \big \{ \varphi^{n_m}_{i_m} + \xi_m -g(y_m),\varphi^{n_m}_{i_m} + \xi_m -H^c_{sup}g(y_m) \big \},\varphi^{n_m}_{i_m} + \xi_m -H^\chi_{inf}g(y_m) \big\} \\ 
    & = \max \big \{ \min \big \{ \varphi^{n_m}_{i_m} + \xi_m -g(y_m),\varphi^{n_m}_{i_m} + \xi_m -(H^c \vec{g})_{i_m} +  O((\Delta x)^2) \big \},\varphi^{n_m}_{i_m} + \xi_m -H^\chi \vec{g} +  O((\Delta x)^2)  \big\}, 
\end{align*}
where we employed the fact that there is $ O((\Delta x)^2) $ error in the intervention operators  approximation due to linear interpolation. Hence, We can without loss of generality assume that $V^h$ is a solution of the scheme so that $V^h(\tau^0,x_i) = g(x_i)$ for all $i$, which is the terminal condition. Hence, we define, 
\begin{equation*}
     S^{(5)}_m =  \max \big \{ \min \big \{ \varphi^{n_m}_{i_m} + \xi_m -g(y_m),\varphi^{n_m}_{i_m} + \xi_m -(H^{c} V^{n_m-1,h_m})_{i_m} +  O((\Delta x)^2) \big \},\varphi^{n_m}_{i_m} + \xi_m - (H^{\chi}V^{n_m-1,h_m})_{i_m} +  O((\Delta x)^2)  \big\} ,
\end{equation*}
and $S^{(4)}_m = S^{(5)}_m$ when $n_m = 0$. Therefore, by (\ref{448}) 
\begin{equation}\label{458}
    S \Big(h_m,(s_m,y_m),\varphi + \xi_m, \mathcal{H}^{h_m}_{sup}V^{h_m},\mathcal{H}^{h_m}_{inf} V^{h_m}\Big) \geq \min \big (  \max\big \{\min[S^{(1)}_m,S^{(2)}_m], S^{(3)}_m\big \}, S^{(5)}_m \big).
\end{equation}
Moreover using Lemma \ref{consistency1} and \ref{consistency2} 
\begin{equation}\label{459}
    \liminf \limits_{m \to \infty}  S^{(5)}_m \geq \max \big \{ \min \big \{ \varphi(t,x) -g(x),\varphi(t,x) -  H_{sup}^c \overline{v}(t,x) \big \},\varphi(t,x) - H_{inf}^{\chi} \overline{v}(t,x) \big\}.
\end{equation}
Taking both  sides of (\ref{458}) to their limit inferiors and applying (\ref{450}), (\ref{451}),(\ref{4511}) and (\ref{459}),
\begin{align}
    \begin{split}
        \liminf \limits_{m\to \infty}  S \Big(&h_m,(s_m,y_m),\varphi + \xi_m , \mathcal{H}^{h_m}_{sup}V^{h_m},\mathcal{H}^{h_m}_{inf} V^{h_m}\Big)\\ 
        &\geq  \min \Bigg( \max \Big\{ \min\big\{ -\varphi_t(t,x)  - b(t,x)\varphi_x(t,x) -  \cfrac{1}{2}\sigma(t,x)^2 \varphi_{xx}(t,x)  - f(t,x) , \varphi(t,x)  - H_{sup}^c \overline{v}(t,x)  \big \} ,\\
    & \qquad  \qquad \varphi(t,x)  - H_{inf}^{\chi} \overline{v}(t,x) \Big\}, 
         \max \big \{ \min \big \{ \varphi(t,x) -g(x),\varphi(t,x) -  H_{sup}^c \overline{v}(t,x) \big \},\varphi(t,x) - H_{inf}^{\chi} \overline{v}(t,x) \big\}  \Bigg ) \\
        & \qquad \qquad \qquad \qquad\qquad \qquad \qquad \qquad \qquad = F_{*}(t,x,\phi(t,x),D\phi(t,x),D^2\phi(t,x),
H^c_{sup} \overline{v}(t,x),H^\chi_{inf} \overline{v}(t,x)),
    \end{split}
\end{align} 
which establish (\ref{454}) as desired. \\ 
It remains to establish (\ref{455}) in the case of $t=T$. Since $s_m \to T$, it is possible that $s_m =T$ (i.e $n_m =0$) for one or more indices $m$ in the sequence. Therefore by (\ref{448}) 
\begin{align}\label{456}
\begin{split}
    S \Big(h_m,(s_m,y_m),\varphi + \xi_m, \mathcal{H}^{h_m}_{sup}V^{h_m},\mathcal{H}^{h_m}_{inf} V^{h_m}\Big) & \leq   \max\big \{\min\{S^{(1)}_m,S^{(2)}_m\}, S^{(3)}_m, S^{(4)}_m \big \} \\
    &= \max\Bigg ( \max\bigg \{\min \big\{S^{(1)}_m,S^{(2)}_m \big\}, S^{(3)}_m\bigg \}, S^{(4)}_m  \Bigg ).\\
    \end{split}
\end{align}
Moreover, using Lemma \ref{consistency1}
and \ref{consistency2} 
\begin{equation}\label{457}
    \limsup \limits_{m \to \infty}  S^{(4)}_m \leq \max \big \{ \min \big \{ \varphi(t,x) -g(x),\varphi(t,x) -  H_{sup}^c \underline{v}(t,x) \big \},\varphi(t,x) - H_{inf}^{\chi} \underline{v}(t,x) \big\}.
\end{equation}
Taking both  sides of (\ref{456}) to their limit superiors and applying (\ref{450}), (\ref{452}), (\ref{4522}) and (\ref{457}),
\begin{align}
    \begin{split}
        \limsup \limits_{m\to \infty}  S \Big(&h_m,(s_m,y_m),\varphi + \xi_m , \mathcal{H}^{h_m}_{sup}V^{h_m},\mathcal{H}^{h_m}_{inf} V^{h_m}\Big)\\ 
        &\leq  \max \Bigg( \max \Big\{ \min\big\{ -\varphi_t(t,x)  - b(t,x)\varphi_x(t,x) -  \cfrac{1}{2}\sigma(t,x)^2 \varphi_{xx}(t,x)  - f(t,x) , \varphi(t,x)  - H_{sup}^c \underline{v}(t,x)  \big \} , \\ &  \qquad  \qquad \varphi(t,x)  - H_{inf}^{\chi} \underline{v}(t,x) \Big\}, 
         \max \big \{ \min \big \{ \varphi(t,x) -g(x),\varphi(t,x) -  H_{sup}^c \underline{v}(t,x) \big \},\varphi(t,x) - H_{inf}^{\chi} \underline{v}(t,x) \big\}  \Bigg ) \\
        & \qquad \qquad \qquad \qquad\qquad \qquad \qquad \qquad \qquad = F^{*}(t,x,\phi(t,x),D\phi(t,x),D^2\phi(t,x),
H^c_{sup} \underline{v}(t,x),H^\chi_{inf} \underline{v}(t,x)),
    \end{split}
\end{align} 
which establishes (\ref{455}) as needed.
\end{proof}
\begin{cor} The solution
$V^{h}$ of (\ref{sch}) converges locally uniformly to the unique  solution of  (\ref{eq:HJBI}), as  $h \downarrow 0$.
\end{cor}

\begin{rem} \label{remarkomicron} 
The idea behind adding the term $\omicron_i^n$ is as follows. We assume that the drift $b(t,x)$ and the forcing term $f(t,x)$ are  given in a way that (i) $b(t,x) = b^{w^{'}}(t,x)$ and (ii) $f(t,x) = f^{w^{'}}(t,x)$ such that  
$$
   \mathcal{L}V + f = \mathcal{L}(w^{'})V + f^{w^{'}}  =   \sup \limits_{w \in W \subset \mathbb{R}}\big(b(t,x,w)\mathcal{D}_x V(t,x) +\cfrac{1}{2}\sigma(t,x)^2\mathcal{D}^2_{xx}V(t,x)+ f(t,x,w)\big ), 
$$
with   $b(t,x,w) = \hat{b}(t,x) + \hat{\hat{b}}(t,x,w)$  and  $f(t,x,w) = \hat{f}(t,x) + \hat{\hat{f}}(t,x,w)$  where the notation $\hat{.}$ and $\hat{\hat{.}}$ are  use for  (sufficiently regular) controlled and uncontrolled components of the functions.\\
Consider the generator $\hat{\mathcal{L}}$ corresponding to an uncontrolled SDE: 
$$
\hat{\mathcal{L}}V(t,x) =   \mathcal{L}(w)V(t,x) - \langle \hat{\hat{b}}(t,x,w), \mathcal{D}_x V(t,x)\rangle.
$$
Letting $X := X(t)$ be a trajectory  satisfying $X(t^n) = x_i$ and $dX(t) = \hat{\hat{b}}(t,X(t),w)dt $ on $(t^n, t^{n+1} ] $   so that $X(t^{n+1}) \approx X(t^n) + \hat{\hat{b}}(t^n,X(t^n),w)\Delta t $.  We define the Lagrangian derivative with respect to $X$ as 
$$
\cfrac{DV}{Dt}(t,X(t),w) := \cfrac{\partial}{\partial t} \big [  V(t,X(t))\big] =   \cfrac{\partial V}{\partial t}(t,X(t)) + \langle \hat{\hat{b}}(t,X(t),w), D_x V(t,X(t))\rangle .
$$
We substitute $\cfrac{DV}{Dt}$ in (\ref{introHJB}) to  get 
$$
\max\Big\{ \min\Big[- \sup \limits_{w \in W \subset \mathbb{R}} \big (\cfrac{DV}{Dt}  +\hat{\mathcal{L}}V+f^w\big ),V-H^c_{sup} V\Big],V-H^\chi_{inf} V \Big\}=0,
$$

which discretization gives
\begin{align*}
\min\Big\{ \max\Big[ \max \limits_{w \in W \subset \overline{\Omega}} \big\{\tilde{V}^{n-1}(t^n, x_i+  \hat{\hat{b}}(t^n,X(t^n),w)\Delta t ) +  &\hat{\mathcal{L}}V^n \Delta t + (\hat{f}^n_i +\hat{\hat{f}}^{n+1}(w))\Delta t \big \}, [H^c_{sup} V^{n-1} -V^n_i]_i + \\ & ([\hat{\mathcal{L}}V^n]_i + \hat{f}^n_i) \Delta t\Big],[H^\chi_{inf} V^{n-1} - V^n_i]_i + ([\hat{\mathcal{L}}V^n]_i + \hat{f}^n_i)\Delta t\Big\}=0,
\end{align*}
and we have $\omicron_i^n = ( [\hat{\mathcal{L}}V^n]_i + \hat{f}^n_i)\Delta t$.\\
As we will see below, adding the term $\omicron_i^n$ allow us to express the scheme as a linear system of
equations at each timestep. However it  introduces
an additional source of discretization error. Then after some simplifications, we can isolate all $V^n$ terms to one side of the equation
\begin{equation} \label{explicit-imp}
    \begin{aligned}[b]
    V^n_i -\cfrac{1}{2}h(\sigma^n_i)^2(\mathcal{D}^2_{xx}V^n)_i = \min \Big\{\max \Big[ ( \tilde{V}^{n-1}(\cdot,x_i + hb^n_i) + hf^n_i )  &, ({H}^{c}V^{n-1})_i  \Big] 
                                                           ,  ({H}^{\chi}V^{n-1})_i  \Big\}.
    \end{aligned}
\end{equation}


Then  we can easily  write (\ref{explicit-imp}) as a system of linear equation (where equation (\ref{explicit-imp}) is the $i$-th row of the linear equation): 
\begin{equation}{\label{sys_lag}}
AV^n = y,
\end{equation}
and $y$ is the $\mathbb{R}^M$ vector whose $i$-th component is the right hand side of (\ref{explicit-imp})
$$
y_i = \min \Big\{\max \Big[ ( \tilde{V}^{n-1}(\cdot,x_i + hb^n_i) + hf^n_i ) , (H^{c}V^{n-1})_i  \Big] 
                                                           ,  (H^{\chi}V^{n-1})_i  \Big\},
$$
and 
\begin{equation}
    (AV^n)_i=  ( V^n_i -\cfrac{1}{2}h(\sigma^n_i)^2(\mathcal{D}^2_{xx}V^n)_i ).
\end{equation}
with $\qquad A=   I -\cfrac{1}{2} h \text{ diag} ((\sigma^n_0)^2,\dots,(\sigma^n_M)^2)\mathcal{D}^2$ given by :
\begin{equation}
    A =  I  + \cfrac{h}{2(\Delta x)^2}
    \begin{pmatrix}
0 & & & &\\
-(\sigma^n_1)^2 & 2(\sigma^n_1)^2 & -(\sigma^n_1)^2   & & & &\\
 & -(\sigma^n_2)^2 & 2(\sigma^n_2)^2 & -(\sigma^n_2)^2  & & &\\
 &  & \ddots &  \ddots& \ddots &   & \\
 &  &  &-(\sigma^n_{M-2})^2  & 2(\sigma^n_{M-2})^2& -(\sigma^n_{M-2})^2 &\\
 &  &  & & -(\sigma^n_{M-1})^2 &2(\sigma^n_{M-1})^2 & -(\sigma^n_{M-1})^2 \\
 &  &  & & & & 0 \\

\end{pmatrix}.
\end{equation}  

Which is a strictly diagonally dominant (Definition \ref{diagdom})  and we can deduce that 
$\mathbb{L}$ is also monotone since  the strict diagonal dominance implies nonsingularity and a monotone matrix is nonsingular [\cite{[AZh]}, Theorem 3.2.5]. As a result the linear system ($\ref{sys_lag}$) has a unique solution
\begin{equation}
    V^n = A^{-1} y.
\end{equation}
\end{rem}

\section{Policy  Iteration }  \label{section3}
\subsection{Policy Iteration framework}
In this section, we consider the scheme (\ref{lagrange11st}) on a  fixed grid 
and study a computational algorithm  to find the approximate value function, which is built on an \textit{approximation in the policy
space} (or \textit{policy iteration}). For generality this approach does not rely on any particular discretization scheme, however we will discuss some general assumptions these scheme should satisfy
for the results presented to hold. 
\no\\
\textbf{Notation.} Henceforward, we work on a discrete finite  space grid
$$
\mathbb{K}: x_0 < x_1 < \cdots < x_{M-1} < x_M  \qquad \text{\textit{Card}}(\mathbb{K})=M+1.
$$
    The \textit{Policy iteration} or \textit{Howard's} algorithm  for discrete variational inequalities  was originally developed
in (\cite{[RBELL]},\cite{[HOWARD]}). The method has been  afterwards extended to variational inequalities in which the obstacle
depends on the solution itself (see \cite{ [ALLA],[AZh],[JPCM]} and the references therein) and to variational inequalities involving multiple obstacle (\cite{[BOKA]}, section 5). Compared to a classical \textit{value iteration} procedure that iterates over value functions, \textit{policy} \textit{iteration} iterates  also over policies themselves, generating  a strictly improved policy in each iteration (except if the iterated policy is already optimal) until convergence. \\
We want to  solve the double obstacle \textit{max-min} problem  of the form
\begin{equation}\label{min-max1}
 \text{Find} \qquad V\in  \mathbb{R}^{\mathbb{K}},\qquad \max\limits_{\beta\in \mathcal{B}^{\mathbb{K}}} \min\limits_{\alpha \in \mathcal{A}^{\mathbb{K}}} \{ A(\alpha,\beta) V- y(\alpha,\beta)\}=0,
\end{equation} 
where $\mathcal{A}$ and $\mathcal{B}$ are non-empty compact set  and for every $\alpha = (\alpha_0, \dots, \alpha_{M} ) \in \mathcal{A}^{\mathbb{K}}$ and $\beta = (\beta_0,\dots , \beta_{M}) \in \mathcal{B}^{\mathbb{K}}$, the matrix $A(\alpha,\beta)\in \mathbb{R}^{{\mathbb{K}} \times {\mathbb{K}}}$ is  monotone, and $y(\alpha, \beta)$ is a vector in $\mathbb{R}^{\mathbb{K}}$, and for 
 which it is understood that : 
 \begin{enumerate}[(i)]
     \item  The controls are row-decoupled, that is, $[A(\alpha,\beta)]_{ij}$ and $[y(\alpha,\beta)]_{i}$ depends only on $(\alpha_i,\beta_i)$.
     \item The order in $\mathbb{R}^{\mathbb{K}}$ and $\mathbb{R}^{\mathbb{K}\times\mathbb{K}}$ is element-wise (i.e : for $x,y$ in $\mathbb{R}^{\mathbb{K}}$, $x\geq y$ if and only if $x_i \geq y_i$ for all $i$).
     \item The $\max$ and $\min$ in (\ref{min-max1}) are with respect to the element-wise order.
 \end{enumerate}
We also assume that
\begin{ass}\label{ass2}
\textbf{(A1)}  The functions $A$ and $y$  are bounded.\\
\textbf{(A2)}  $(\alpha,\beta) \to A(\alpha,\beta)^{-1}$ is nonsingular and monotone.\\
\textbf{(A3)} For each $(\alpha,\beta)\in \mathcal{A}^{\mathbb{K}} \times \mathcal{B}^{\mathbb{K}}$, $A(\alpha,\beta)$ is a strictly diagonally dominant  M-matrix  with non-negative diagonals.

\end{ass}
\begin{axiom}\label{diagdom} We say  row $i$ of a complex  matrix $A:=(a_{ij})$ is weakly diagonally dominant  if $|a_{ii}| \geq \sum \limits_{i \neq j} |a_{ij}|$.\\
 We say row $i$ of a complex matrix $A:=(a_{ij})$ is strictly diagonally dominant if $|a_{ii}| > \sum \limits_{i \neq j} |a_{ij}|$.\\ 
We say $A$ is strictly diagonally dominant if all its rows are SDD. The same holds for WDD. \\
A complex matrix $A$ is said to weakly chained diagonally dominant (WCDD) if $A$ is  WDD, and for each row $i$  there exists a  path in the graph of $A$ from $i$ to an SDD row $i^{'}$.\\
$A$ is said to be monotone if and only if $A$ is invertible and $A^{-1} \geq 0$ componentwise. Moreover $A$ is nonsingular if $AA^{-1} =A^{-1}A= I$ where $I$ is the identity matrix.
\end{axiom}
We recall  some results used further in the paper for the reader convenience.
\begin{theo} \label{theomono}  (\cite{[AZh]}, Theorem 3.2.5.) If $A$ is a square w.d.d. Z-matrix with nonnegative diagonals, then the
following are equivalent:
    \begin{enumerate}[(i)]
    \item $A$ is a w.c.d.d. matrix.
        \item $A$ is an M-matrix.
        \item $A$ is a nonsingular matrix.
    \end{enumerate}
\end{theo}

For problems such as (\ref{min-max1}), we introduce the general  policy iteration algorithm below.
\begin{enumerate}
    \item \textbf{Initialization}
    \begin{enumerate}
        \item \textbf{Define Sets}:
        \begin{itemize}
            \item $\mathcal{A}^{\mathbb{K}}$ the set of all possible $\alpha$ policies and $\mathcal{B}^{\mathbb{K}}$ the set of all possible $\beta$ policies.
        \end{itemize}
        \item \textbf{Initial Policy}:
        \begin{itemize}
            \item Select an initial policy $(\alpha^0, \beta^0) \in \mathcal{A}^{\mathbb{K}} \times \mathcal{B}^{\mathbb{K}}$.
        \end{itemize}
        \item \textbf{Initial Value Vector}:
        \begin{itemize}
            \item Initialize $V$ as a vector in $\mathbb{R}^{\mathbb{K}}$, e.g., $V^0 = 0$.
        \end{itemize}
    \end{enumerate}

    \item \textbf{Policy Evaluation}

    Given a policy $(\alpha^k, \beta^k)$, solve for $V^{k}$ in the  system:
    \begin{equation*}
    A(\alpha^k, \beta^k) V^{k} = y(\alpha^k, \beta^k).
    \end{equation*}
   
    \item \textbf{Policy Improvement}

    Update the policies $\alpha$ and $\beta$ by solving the following optimization problems:
    \begin{enumerate}
        \item \textbf{Update $\alpha$ Policy}:
        \begin{equation}
        \alpha^{k+1} = \arg \min_{\alpha \in \mathcal{A}^{\mathbb{K}}} \{ A(\alpha, \beta^k) v^{k} - y(\alpha, \beta^k) \}.
        \end{equation}
        \item \textbf{Update $\beta$ Policy}:
        \begin{equation}
        \beta^{k+1} = \arg \max_{\beta \in \mathcal{B}^{\mathbb{K}}} \{ A(\alpha^{k+1}, \beta) v^{k} - y(\alpha^{k+1}, \beta) \}.
        \end{equation}
    \end{enumerate}

    \item \textbf{Convergence Check}

    Check if the policies $\alpha$ and $\beta$ have converged:
    \begin{equation}
    \max_{\beta \in \mathcal{B}^{\mathbb{K}}} \min_{\alpha \in \mathcal{A}^{\mathbb{K}}} \{ A(\alpha, \beta) v^{k} - y(\alpha, \beta) \} \leq \epsilon,
    \end{equation}
    where $\epsilon$ is a small positive tolerance.

    \item \textbf{Iteration}\\
    If the convergence has not been achieved, set $k = k + 1$ and return to Step 2.
\end{enumerate}

\subsection{Application to the HJBQVI}\label{applicationij}

Recall the discretized numerical scheme $S$ introduced in (Proposition \ref{proposition31}):
\
\begin{align}\label{policyeq}
\begin{split}
\max\limits_{p \in \{0,1\}} 
       \Bigg\{  \min\limits_{q\in \{ 0,1\}} \bigg[-(1-p)(1-q)(\tilde{V}^{n-1}(\cdot,x_i + hb^n_i) -& V^n(\cdot,x_i) +\cfrac{1}{2}h(\sigma^n_i)^2(\mathcal{D}^2V^n)_i+hf^n_i )-\\&
       (1-p)q ((H^{c}V^{n-1})_i -V^n_i+ \omicron_i^n ) \bigg] - 
        p ((H^{\chi}V^{n-1})_i -V^n_i + \omicron_i^n  ) \Bigg\} =0.
\end{split}
\end{align}
To write this scheme in the form of (\ref{min-max1}),
we define the matrix-valued function $A$ and vector-valued function $y$ such that for each
vector $V = (V_0, \cdots , V_M)$ and integer $i$ satisfying $0\leq i \leq M$,
\begin{align*}
[A(\alpha,\beta)V]_i = \big ( V_i - \cfrac{1}{2}h(\sigma^n_i)^2(\mathcal{D}^2V)_i \big ) ,
\end{align*}
and 
\begin{align*}
[y(\alpha,\beta)]_i = (1-p)(1-q)\big (\tilde{V}^{n-1}(\cdot,x_i + hb^n_i) +hf^n_i \big ) + (1-p)q ((H^{c}V^{n-1})_i )+ p ((H^{\chi}V^{n-1})_i) ,
\end{align*}
with 
$$
\alpha_i = (q_i,\xi_i) , \qquad \text{and} \qquad  \beta_i = (p_i,\eta_i), \qquad (p,q) =\{0,1\}\times\{0,1\}.
$$
\begin{lem} $A(\alpha,\beta)$ is a monotone matrix.
\end{lem}
\begin{proof} We will prove that $A(\alpha,\beta)$ is an s.d.d. Z-matrix with positive diagonals so that by so that  by Theorem \ref{theomono} it is an M-matrix.
\begin{align}\label{offdiag}
[A(\alpha,\beta)V]_i = ( 1 + \dfrac{h(\sigma^n_i)^2}{(\Delta x)^2} ) V_i - \dfrac{1}{2}\dfrac{h(\sigma^n_i)^2}{(\Delta x)^2}V_{i+1} - \dfrac{1}{2}\dfrac{h(\sigma^n_i)^2}{(\Delta x)^2}V_{i-1} .
\end{align}
From (\ref{offdiag}), $A(\alpha,\beta)$ ) has nonpositive off-diagonals and as such, is a Z-matrix. Now setting $V= \Vec{e}$, for all $i$, 
$$
\sum\limits_{j} [A(\alpha,\beta)V]_{ij} = [A(\alpha,\beta)\Vec{e}] = 1 > 0.
$$
Moreover, $$[A(\alpha,\beta)]_{ii} =  ( 1 + \dfrac{h(\sigma^n_i)^2}{(\Delta x)^2} ) > 1,$$
from  which it follows that $A(\alpha,\beta)$  is s.d.d. with positive diagonals since 
$$
[A(\alpha,\beta)]_{ii} > - \sum\limits_{j\neq i} [A(\alpha,\beta)V]_{ij} =  \sum\limits_{j\neq i} |[A(\alpha,\beta)V]_{ij}| \geq 0.
$$
\end{proof}
\subsection{Computational algorithm }
In this part  we give two computational algorithms  to find the approximate value function $V^h$ of the scheme $S$ from which we  deduce the optimal controls for  each player (the NE strategy) as
well as the optimal evolution of the state. The algorithm \ref{alg1} describes the implementation of the optimal impulses controls to perform the desired computations based on\textit{ Policy iteration}. 

\textbf{Part 1}  Computes the  values in the whole time–space grid
\begin{itemize}
    \item Step 1:  \textbf{Recursive computation}. From  the terminal value $V^0 = V^h(T , x) = g(x)$,
we recursively compute $V^h(t^n, x_i)$ for $(t^n, x_i)$ in a given time–space grid when no player intervenes by taking $p=q=0$:
\begin{align*}
\begin{split}
\left \{
\begin{array}{l}
[A(\alpha,\beta)V]_i = \big ( V_i - \cfrac{1}{2}h(\sigma^n_i)^2(\mathcal{D}^2V)_i \big ) ,\\ \\
 y(\alpha,\beta)_i =  \tilde{V}^{n-1}(\cdot,x_i + hb^n_i) +hf^n_i .
\end{array}
\right .
\end{split}
\end{align*}
This requires that for any time–space grid point $(t^n, x_i)$ the quantity $(x_i + hb^n_i)$ remains in the space domain. We use a linear interpolation to  compute $V^h$ at such space points. 
\item  Step 2 : \textbf{The Minimizing player-$\eta$ Intervenes}. This player  intervenes with at time $t^n$ only when 
$$
H^{\chi}V(t^{n-1}, x^{*}_i) + \cfrac{1}{2}h(\sigma^n_i)^2\mathcal{D}^2V(t^n, x^{*}_i)  <  V(t^n, x^{*}_i),
$$
with the optimal state $x^{*}$. The new approximate value $V^h$ is computed by  taking $p=1$.
\begin{align*}
\begin{split}
\left \{
\begin{array}{l}
[A(\alpha,\beta)V]_i = \big ( V_i - \cfrac{1}{2}h(\sigma^n_i)^2(\mathcal{D}^2V)_i \big ) ,\\ \\
 y(\alpha,\beta)_i =    (H^{\chi}V^{n-1})_i =  \inf \limits_{
 \eta\in \mathbb{K}}[\tilde{V}(\tau^{n-1},x+\eta) +\chi(\tau^n,\eta)].
\end{array}
\right .
\end{split}
\end{align*}
For the above system a
linear interpolation is used to  compute   $V^h$ at $x+\eta$. The policy iteration of Algorithm \ref{alg1} will give the optimal size of the impulse controls.
\item  Step 3: \textbf{The maximizing player-$\xi$ Intervenes}. When the minimizing player does not intervene at time  $t^n$, the
maximizing player  might intervene when  
$$
H^{c}V(t^{n-1}, x^{*}_i) + \cfrac{1}{2}h(\sigma^n_i)^2\mathcal{D}^2V(t^n, x^{*}_i)  > V(t^n, x^{*}_i), 
$$ with the optimal state $x^{*}$ and the new approximate value $V^h$ has to be  computed by  taking $p=0$ and $q=1$.
\begin{align*}
\begin{split}
\left \{
\begin{array}{l}
[A(\alpha,\beta)V]_i = \big ( V_i - \cfrac{1}{2}h(\sigma^n_i)^2(\mathcal{D}^2V)_i \big ) ,\\ \\
 y(\alpha,\beta)_i =     (H^{c}V^{n-1})_i =  \sup \limits_{ \xi \in \mathbb{K}}[\tilde{V}(\tau^{n-1}, x_i + \xi)-c(\tau^n, \xi)].
\end{array}
\right .
\end{split}
\end{align*}
\end{itemize}
\textbf{Part 2:}  Computes the optimal impulse controls strategy pair  $\{\psi = (u^* = (\tau^*_m,\xi^*_{m}),v^*=(\rho^*_l,\eta^*_{l})\}$  ,  the optimal state $x^*$  evolution  and the approximate value function  $V^h$.

\begin{algorithm}[H]
\caption{Policy Iteration Algorithm for Scheme $S$} \label{alg1}
\begin{algorithmic}[1]
\State \textbf{Input:} 
\State Time–space grid points $(t_n,x_i)$ , where $t_n \in \{ t_0  = T, t_1, . . . , t_N = 0 \}$ and $x_i \in  \{x_0, x_1, . . . , x_k, . . . , x_M\}$  ;
\State Maximum iterations $N$, discretization step $h$;
\State Terminal value $g(x)$ ;
\State Functions $\sigma$, $b$,  $f$;
\State Initial policies $\alpha^0,\beta^0$;
\State Initial impulse policies $\xi^0 = [\xi^0_1, \cdots, \xi^0_M]$ and $\eta^0 = [\eta^0_1,\cdots,\eta^0_M]$
\State Functions $c,\chi$;
\State Set tolerance $\epsilon$;
\State \textbf{Output:} Value function $V$,Optimal policies $\alpha$, $\beta$, optimal state evolution $x^*$,optimal controls  $(\tau^*_m,\xi^*_{m})$ and $(\rho^*_l,\eta^*_{l})$;
\State Initialize $V^0_i \leftarrow g(x_i)$ for  all $i = 0, \cdots , M $;
\State Initialize policies $\alpha^0 \leftarrow 0$, $\beta^0 \leftarrow 0$ for all $i$ in the space grid
\State Initialize  impulses values $\tau^*_0,\xi^*_0$, $\rho^*_0,\eta^*_0 \leftarrow 0,0,0,0$ \\
\State \textbf{Part 1}
\State \textbf{for} $k = 1$ to $N$ \textbf{do}
 \State \-\hspace{0.5cm} $V_{\text{prev}} \leftarrow V^{k-1}$
    \State \-\hspace{0.5cm}  \textbf{Policy Evaluation:}
    \State \-\hspace{0.5cm} \textbf{for} $i=0$  to $M$ do
     \State \-\hspace{1.3cm} Define $A(\alpha^{k-1}_i,\beta^{k-1}_i) V_i \leftarrow V(x_i) - \frac{1}{2}h(\sigma_i)^2 (\mathcal{D}^2_{xx} V)_i$
     \State \-\hspace{1.3cm} \textbf{Step 1}
        \State \-\hspace{1.3cm}  Compute $\text{term1}^k_i$: $V_{\text{prev}}(x_i + hb)  + h f_i$
        \State \-\hspace{1.3cm} \textbf{Step 2} \Comment{Policy evaluation and improvement over impulse policies $\eta^k$}
       \State \-\hspace{1.3cm}  Compute $\text{term2}^k_i$: $ (H^{\chi}V_{\text{prev}})_i$
       \State \-\hspace{2cm} Find $\eta^k_i \leftarrow  \arg \inf \limits_{
 \eta\in \mathbb{K}}[\tilde{V}_{\text{prev}}(x+\eta) +\chi(t^k,\eta)].$
       \State \-\hspace{2cm} $\text{term2}^k_i \leftarrow [\tilde{V}_{\text{prev}}(x+\eta^k_i) +\chi(t^k,\eta^k_i)] $ 
        \State \-\hspace{1.3cm} \textbf{ Step 3} \Comment{Policy evaluation and improvement over impulse policies $\xi^k$}
       \State \-\hspace{1.3cm}  Compute $\text{term3}^k_i \leftarrow (H^c V_{\text{prev}})_i$ 
       \State \-\hspace{2cm} Find $\xi^k_i \leftarrow  \arg \sup \limits_{
 \xi\in \mathbb{K}}[\tilde{V}_{\text{prev}}(x+\xi) -c(t^n,\xi)].$
       \State \-\hspace{2cm} $\text{term3}^k_i \leftarrow [\tilde{V}_{\text{prev}}(x+\xi^k_i) -c(t^k,\xi^k_i)] $ 
       \State \-\hspace{1.3cm} Define  $y(\alpha^{k-1}_i,\beta^{k-1}_i)_i \leftarrow \min\left( \max(\text{term1}^k_i, \text{term3}^k_i), \text{term2}^k_i \right)$   
        \State \-\hspace{1.3cm} Update $V_i \leftarrow $ solution of $A(\alpha_k,\beta_k) V_i = y(\alpha_k,\beta_k)_i  $
\State \-\hspace{0.5cm} \textbf{end for}
    
    \textbf{if} {$\max |V_k - V_{\text{prev}}| < \epsilon$} \textbf{then} 
        \State \-\hspace{1.3cm}\textbf{Converged}
        \State \-\hspace{1.3cm} \textbf{Break}
    
    \State \-\hspace{0.5cm} \textbf{Policy Improvement:}
    \State \-\hspace{0.5cm} \textbf{for} $i=0$  to $M$ do
         \State \-\hspace{1cm} Update $\alpha^k_i \leftarrow \arg \min_{\alpha} \left[  V(x_i) - \frac{1}{2}h(\sigma_i)^2 (\mathcal{D}^2_{xx} V)_i - (1-\alpha) \big ((V_{\text{prev}}(x_i + h b_i) + h f_i\big ) + \alpha \big ((H^{c}V_{\text{prev}})_i \big ) \right]$
         \State \-\hspace{1cm} Update $\beta^k_i \leftarrow \arg \max_{\beta} \left[ V(x_i) - \frac{1}{2}h(\sigma_i)^2 (\mathcal{D}^2_{xx} V)_i - (1-\beta)(1-\alpha^k)\big ((V_{\text{prev}}(x_i + h b_i) + h f_i\big ) \right.$ 
         \State \-\hspace{6cm} $ + (1-\beta)\alpha^k((H^{c}V_{\text{prev}})_i )+ \beta((H^\chi V_{\text{prev}})_i) \left. \right]$
    \State \textbf{Return} $V, \alpha, \beta$,$\eta$,$\xi$
     \State  \textbf{Part 2}
     \State Initial state $x_1^* \leftarrow x$
     \State $l,m \leftarrow 1,1$
     \State $V^k$ is such that  $V_k \approx V(t^k, x_k^*)$
     \State  \textbf{for} $k=1$  to $N$ \textbf{do}
    \State \-\hspace{0.5cm} \textbf{if} $V^k > H^{\chi}V_{\text{prev}}( x^{*}_k) + \cfrac{1}{2}h(\sigma^n_i)^2\mathcal{D}^2V^k( x^{*}_k)  $ \textbf{then} 
     \State \-\hspace{1.3cm} $\rho^*_l,\eta^*_l \leftarrow t_k, \eta_k $
     \State \-\hspace{1.3cm} $x_{k+1}^* \leftarrow x^*_k + \eta_k$ 
     \State \-\hspace{1.3cm} $V_k^* \leftarrow V(x^{*}_k)$
     \State \-\hspace{1.3cm}  $l \leftarrow l+1$
      \State \-\hspace{0.5cm} \textbf{else if}  $ V^k < H^{c}V_{\text{prev}}( x^{*}_k) + \cfrac{1}{2}h(\sigma^n_i)^2\mathcal{D}^2V(t^k, x^{*}_k)$
    \textbf{then} 
     \State \-\hspace{1.3cm} $\tau^*_m,\xi^*_m \leftarrow t_k, \xi_k $
     \State \-\hspace{1.3cm} $x_{k+1}^* \leftarrow x^*_k + \xi_k$ 
     \State \-\hspace{1.3cm} $V_k^* \leftarrow V(x^{*}_k)$
     \State \-\hspace{1.3cm}  $m \leftarrow m+1$
     \State \-\hspace{0.5cm} \textbf{else} 
     
     \State \-\hspace{1.3cm} $x_{k+1}^* \leftarrow x^*_k + hb(t_k, x^*_k)$ 
     \State \-\hspace{1.3cm} $V_k^* \leftarrow V(x^{*}_k)$
    \State \-\hspace{0.5cm} \textbf{end if} 
    \State \textbf{end for}
\end{algorithmic}
\end{algorithm}

\section{Implementation and computational  studies}\label{section5}
\no
\subsection{Computing data for the exchange rate problem}
For the numerical implementation, we  use the parameters in table \ref{table1}. 
\begin{table}[h!]
\centering
\begin{tabular}{p{2cm}  p{3cm} p{4cm} p{2.5cm}}
\hline 
    Parameter&  Desc. &  Game 1 ($x_0 = 2.5$)  \\
    \hline 
    t & initial time & 0 \\
    T &  Horizon &  1  \\ 
    h & step & 0.05   \\ 
    \hline 
   $\mu$ & Drift speed &  0.25  \\
  $\sigma$ &  Volatility &   0.30 \\
    $x^*$ & Target rate & 1    \\ 
    $k_i$ & Fixed cost &  0.1  \\
  $\lambda_i$ &  Proportional cost  &  1\\
    \hline 
  $g(x)$ &  Terminal value &  0  \\
  $f(x)$ & Running reward & $-(x - x^{*})^2$ \\ 
    \hline
 \end{tabular}
 \caption{ Optimal control of the FEX rate: parameters}\label{table1}
 \end{table}
\subsection{Numerical results}
\begin{figure}
\centering
\begin{subfigure}{0.455\linewidth}
    \includegraphics[width=\linewidth]{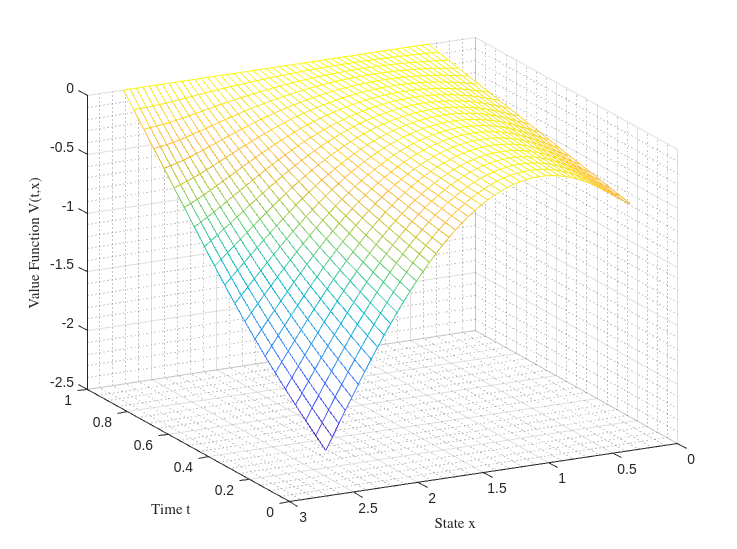}
\end{subfigure}
\begin{subfigure}{0.455\linewidth}
    \includegraphics[width=\linewidth]{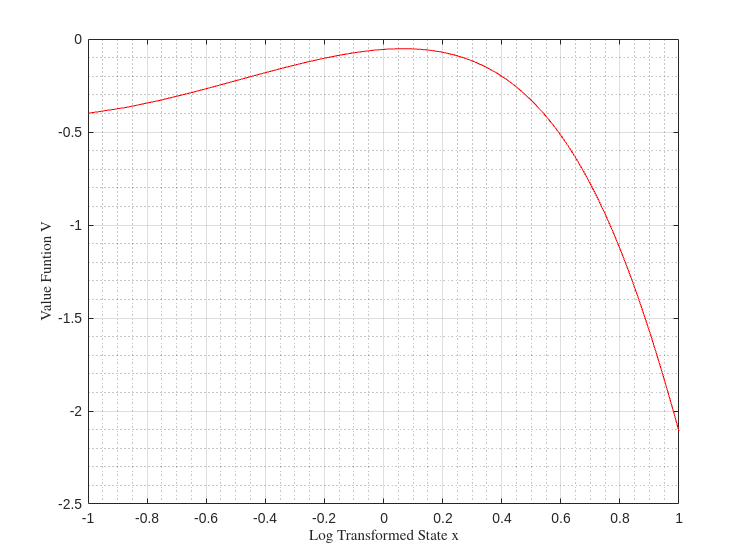}
\end{subfigure}
\caption{Value Function $V(t,x)$ (left) an $V(t=0,x)$ (right) at $t=0$.\label{fig1}}
\end{figure}

\begin{figure}
    \centering
    \includegraphics[width=1\linewidth]{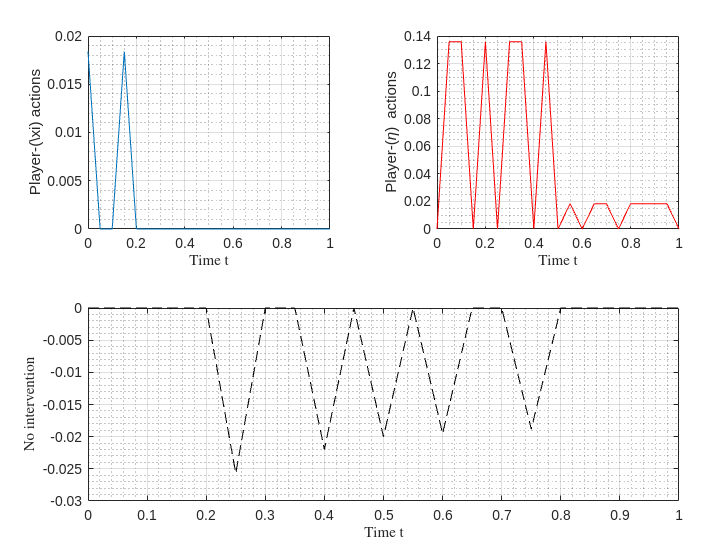}
    \caption{Impulses actions for Player-$\xi$(top-left) and Player-$\eta$ (top-right), and game evolution (bottom) with no impulses.\label{fig2}}
    \label{fig:enter-label}
\end{figure}

\begin{figure}
\centering
\begin{subfigure}{0.455\linewidth}
    \includegraphics[width=\linewidth]{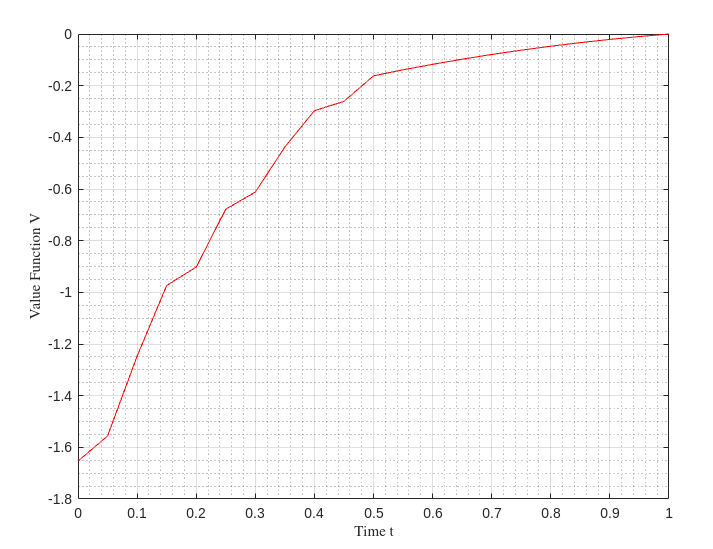}
\end{subfigure}
\begin{subfigure}{0.455\linewidth}
    \includegraphics[width=\linewidth]{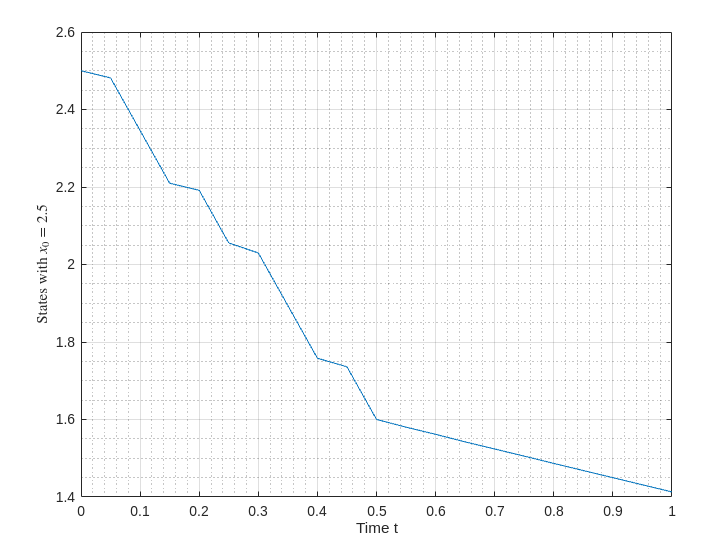}
\end{subfigure}
\caption{Optimal Value Function for the initial state $x_0 = 2.5$ (left) and optimal states starting from $x_0 =2.5$(right). \label{fig3}}
\end{figure}

See  Fig.\ref{fig1}, Fig. \ref{fig2} and Fig. \ref{fig3}.

\section{Conclusion}\label{section6}

In this paper, we have considered a stochastic finite-time horizon, two-player, zero-sum DGs, where both players use impulse controls. The aims were to optimize the discounted terminal gain/cost functional, approximate the value function, and describe an optimal impulses strategies for the two players. After analyzing the related HJBI double-obstacle equation (\ref{eq:HJBI})in the viscosity solution framework provided that it verifies some underlined assumptions and a verification theorem, we
have proposed a discrete-time approximation scheme for this class of DGs given by the approximate equation (\ref{lagrange1st}). We have further demonstrated that the scheme verify some monotonicity, stability and non-local consistency and therefore converges to its initial HJBI as the discretization step $h$ goes to zero.
In general, our numerical approximation involved discretizing the state and time variables and approximating the
derivatives using finite difference. The accuracy of that method depends on the size of the grid and time steps used in the
discretization, as well as the smoothness of the underlying functions. In some cases, it may be necessary to refine the grid and
time steps to achieve a desired level of accuracy. Overall, this numerical approximation of the zero-sum stochastic differential impulse
control games requires a careful balance between accuracy and computational efficiency, as the complexity of the problem can
quickly escalate with the size of the state space and the duration of the game.\\
Our major contributions are  the convergence result
for the approximate value function, and a policy iteration based algorithm to  find the value function in a practical implementation. Moreover, we have given some meaningful dynamics $X_t, b(\cdot),\sigma(\cdot), f(\cdot), c(\cdot),\chi(\cdot)$ to
 apply our results to a foreign exchange rate (FEX) problem  game where a central bank and a commercial institution intervene in the market.
In such situation our results have been successfully implemented to derive the value functions as well as the impulses strategies for both actors. 
We intend to develop this work in  more  directions in the future:
\begin{enumerate}
    \item  Consider the problem  in an infinite time horizon with no backward computation and no terminal value function. 
    \item Adopt  a direct control discretization scheme for solving the HJBI. This will  require  revisiting the monotonicity assumption of the matrices $A(\alpha,\beta)$ for singularity in the policy iteration framework.
\end{enumerate}

\end{document}